\documentclass[10pt]{article}
\usepackage{epsf}
\usepackage{epsfig}

\def\PP{\rm \hbox{I\kern-.2em\hbox{P}}}
\def\RR{\rm \hbox{I\kern-.2em\hbox{R}}}
\def\NN{\rm \hbox{I\kern-.2em\hbox{N}}}
\def\ZZ{\rm {{\rm Z}\kern-.28em{\rm Z}}}
\def\CC{\rm \hbox{C\kern -.5em {\raise .32ex \hbox{$\scriptscriptstyle
|$}}\kern-.22em{\raise .6ex \hbox{$\scriptscriptstyle |$}}\kern .4em}}
\def\vp{\varphi}
\def\<{\langle}
\def\>{\rangle}
\def\ti{\tilde}

\def\e{\varepsilon}
\def\sm{\setminus}
\def\nl{\newline}

\def\bq{{\bf q}}
\def\cT{{\cal T}}

\def\cD{{\cal D}}

\def\cP{{\cal P}}

\def\cO{{\cal O}}
\def\cH{{\cal H}}
\def\bn{{\rm \bf n}}
\def\bt{{\rm \bf t}}
\def\Chi{\raise .3ex\hbox{\large $\chi$}} \def\vp{\varphi}
\def\lsima{\hbox{\kern -.6em\raisebox{-1ex}{$~\stackrel{\textstyle<}{\sim}~$}}\kern -.4em}
\def\lsim{\hbox{\kern -.2em\raisebox{-1ex}{$~\stackrel{\textstyle<}{\sim}~$}}\kern -.2em}
\def\({\Bigl (}
\def\){\Bigr )}
\def\({\Bigl (}
\def\){\Bigr )}

\newcommand{\be}{\begin{equation}}
\newcommand{\ee}{\end{equation}}
\newcommand{\bea}{$$ \begin{array}{lll}}
\newcommand{\eea}{\end{array} $$}
\newcommand{\bi}{\begin{itemize}}
\newcommand{\ei}{\end{itemize}}
\newcommand{\iref}[1]{(\ref{#1})}
\newtheorem{theorem}{Theorem}[section]
\newtheorem{remark}[theorem]{Remark}

\newtheorem{proposition}[theorem]{Proposition}

\newtheorem{definition}[theorem]{Definition}

\newtheorem{conjecture}[theorem]{Conjecture}

\def\proof{{\noindent \bf Proof: }}

\def\sq{\hfill $\diamond$\\}
\def\cO{\mathcal O}
\def\R{\mathbb R}

\def\P{{\rm \hbox{I\kern-.2em\hbox{P}}}}
\def\H{{\rm \hbox{I\kern-.2em\hbox{H}}}}

\def\ssdelta{{\scriptscriptstyle\delta}}

\def\n{{\rm \bf n}}
\def\t{{\rm \bf t}}

\def\trans{{\mathrm T}}

\newcommand\sep{\; ; \;}

\def\gsim{\hbox{\kern -.2em\raisebox{-1ex}{$~\stackrel{\textstyle>}{\sim}~$}}\kern -.2em}

\usepackage{amsmath}
\usepackage{amsfonts}
\usepackage{amssymb}

\usepackage{enumerate}

\usepackage{graphicx}
\DeclareMathOperator\Id{Id}

\DeclareMathOperator\tr{tr}
\DeclareMathOperator\proj{P}
\DeclareMathOperator\move{U}

\DeclareMathOperator\TV{TV}

\DeclareMathOperator\disc{disc}

\def\II{\mathrm {I\kern-0.1exI}}

\begin{document}
\title{\bf Anisotropic smoothness classes:\\ \bf from finite element approximation to image models} 
\author{Jean-Marie Mirebeau and Albert Cohen}
\maketitle
\date{}
\begin{abstract}
We propose and study quantitative measures of smoothness 
$f\mapsto A(f)$ which are adapted to anisotropic features such as
edges in images or shocks in PDE's. These quantities
govern the rate of approximation by adaptive finite elements,
when no constraint is imposed on the aspect ratio of the triangles,
the simplest example being $A_p(f)=\|\sqrt{|{\rm det}(d^2f)|}\|_{L^\tau}$
which appears when approximating in the $L^p$ norm by piecewise linear elements
when $\frac 1 \tau=\frac 1 p+1$. The quantities $A(f)$ are not
semi-norms, and therefore cannot be used to define linear function spaces.
We show that these quantities can be well defined by mollification when $f$ has
jump discontinuities along piecewise smooth curves. This motivates
for using them in image processing as an alternative to the frequently
used total variation semi-norm which does not account for 
the smoothness of the edges. 
\end{abstract}

\section{Introduction}

There exists various ways of measuring the smoothness of functions
on a domain $\Omega\subset \RR^d$, 
generally through the definition of an appropriate {\it smoothness space}. 
Classical instances are Sobolev, H\"older and Besov spaces. Such spaces are of common use
when describing the regularity of solutions to partial differential equations.
From a numerical perspective, they are also useful in order to sharply characterize at which
rate a function $f$ may be approximated by simpler functions such as Fourier series,
finite elements, splines or wavelets (see \cite{Co,DL,De} for surveys on
such results).

Functions arising in concrete applications may have
inhomogeneous smoothness properties, in the sense that they exhibit
area of smoothness separated by localized discontinuities. Two typical
instances are (i) edge in functions representing real images and 
(ii) shock profiles in solutions to non-linear hyperbolic PDE's.
The smoothness space that is best taylored to take such features
into account is the space $BV(\Omega)$ of bounded variation functions. This space
consists of those $f$ in $L^1(\Omega)$ such that $\nabla f$ is a bounded measure,
i.e. such that their total variation
$$
\TV(f)=|f|_{BV}:=\max \left\{\int_\Omega f {\rm div}(\vp)\; ; \; \vp\in\cD(\Omega)^d,\; \|\vp\|_{L^\infty}\leq 1\right\}
$$
is finite. Functions of bounded variation are allowed to have jump discontinuities
along hypersurfaces of finite measure. In particular, the characteristic function of a smooth
subdomain $D\subset \Omega$ has finite total variation equal to the $d-1$-dimensional
Hausdorff measure of its boundary:
\be
|\Chi_D|_{BV}=\cH_{d-1}(\partial D).
\label{haus}
\ee
It is well known that $BV$ is a regularity space for certain 
hyperbolic conservation laws \cite{GR,Le}, in the sense that the total variation
of their solutions remains finite for all time $t>0$. This space also plays an
important role in image processing since the seminal paper \cite{FOR}.
Here, a small total variation is used as a prior to describe the mathematical 
properties of ``plausible images'', when trying to restore an unknown
image $f$ from an observation $h=Tf+e$ where $T$ is a known operator
and $e$ a measurement noise of norm $\|e\|_{L^2}\leq \e$. The restored image
is then defined as the solution to the minimization problem
\be
\min_{g\in BV} \{|g|_{BV}\; ; \; \|Tg-h\|_{L^2}\leq \e\}.
\label{minBV}
\ee
From the point of view of approximation theory, it was shown in \cite{CDPX} that
the space $BV$ is almost characterized by expansions in wavelet bases. For example,
in dimension $d=2$, if $f=\sum d_\lambda\psi_\lambda$ 
is an expansion in a tensor-product $L^2$-orthonormal wavelet basis, one has
$$
(d_\lambda)\in \ell^1 \Rightarrow f\in BV \Rightarrow (d_\lambda)\in w\ell^1,
$$
where $w\ell^1$ is the space of weakly summable sequences. The fact that 
the wavelet coefficients of a $BV$ function are weakly summable implies
the convergence estimate
\be
\|f-f_N\|_{L^2}\leq CN^{-1/2}|f|_{BV},
\label{wavN}
\ee
where $f_N$ is the {\it nonlinear} approximation of $f$ obtained by retaining the
$N$ largest coefficients in its wavelet expansion. Such approximation results 
have been further used in order to justify the performance of compression or
denoising algorithms based on wavelet thresholding \cite{CDDD1}.

In recent years, it has been observed that the space $BV$ (and more generally
classical smoothness spaces) do not provide a fully satisfactory description of
piecewise smooth functions arising in the above mentioned applications.
Indeed, formula \iref{haus} reveals that the total variation only takes into account
the {\it size} of the sets of discontinuities and not their geometric {\it smoothness}.
In image processing, this means that the set of bounded variation images does
not make the distinction between smooth and non-smooth edges as long as they have finite length.

The fact that edges have some geometric smoothness can be exploited in order to study
approximation procedures which outperform wavelet thresholding in terms of 
convergence rates. For instance, it is easy to prove that if $f=\Chi_D$ 
where $D$ is a bidimensional domain with smooth boundary, one can find a sequence
of triangulations $\cT_N$ with $N$ triangles such that the convergence estimate
\be
\|f-I_{\cT_N}f\|_{L^2} \leq CN^{-1},
\label{triN}
\ee
holds, where $I_{\cT}$ denotes the 
piecewise linear interpolation operator on a triangulation $\cT$.
Other methods are based on thresholding a decomposition of the 
function in bases or frames which differ from classical wavelets, see e.g. \cite{CD,LM,ACDDM}.
These methods also yield improvements over \iref{wavN} similar to \iref{triN}. The common feature in
all these approaches is that they achieve 
{\it anisotropic refinement} near the edges. For example,
in order to obtain the estimate \iref{triN}, 
the triangulation $\cT_N$ should include a thin layer of triangles
which approximates the boundary $\partial D$. These triangles typically have
size $N^{-2}$ in the normal direction to $\partial D$ and $N^{-1}$ in the tangential direction,
and are therefore highly anisotropic.  

Intuitively, these methods are well adapted to
functions which have anisotropic smoothness properties
in the sense that their local variation is significantly stronger
in one direction. Such properties are not well described by
classical smoothness spaces such as $BV$, and a natural
question to ask is therefore:
\nl
\nl
{\it What type of smoothness properties 
govern the convergence rate of anisotropic refinement
methods and how can one quantify these properties ?}
\nl
\nl
The goal of this paper is to answer this question,
by proposing and studying measures of smoothness which are suggested
by recent results on anisotropic finite element approximation
\cite{BBLS, CSX}. Before going further, 
let us mention several existing approaches which have been developed
for describing and quantifying anisotropic smoothness,
and explain their limitations.
\begin{enumerate}
\item
The so-called {\it mixed smoothness} classes have been introduced 
and studied in order to
describe functions which have a different order of
smoothness in each coordinate, see e.g. \cite{Ni,Te}. These spaces
are therefore not adapted to our present goal since
the anisotropic smoothness that we want to describe 
may have preferred directions that are not 
aligned with the coordinate axes and that may vary from
one point to another (for example an image with a curved edge).
\item
Anisotropic smoothness spaces
with more general and locally varying directions have
been investigated in \cite{KP}. Yet, in such spaces
the amount of smoothness in different directions
at each point is still fixed in advance and therefore
again not adapted to our goal, since this amount
may differ from one function to another (for example
two images with edges located at different positions).
\item
A  class of functions which is often used to 
study the convergence properties of anisotropic approximation
methods is the family of $C^m-C^n$ {\it cartoon images}, i.e.
functions which are $C^m$ smooth on a finite number
of subdomains $(\Omega_i)_{i=1,\cdots,k}$ separated by a union of
discontinuity curves 
$(\Gamma_{j})_{j=1,\cdots,l}$
that are $C^n$ smooth. The defects of this class are revealed
when searching for simple expression that quantifies the amount of
smoothness in this sense. A natural choice is
to take the supremum of all $C^m(\Omega_i)$ norms of $f$ and
$C^n$ norms of the normal parametrization of $\Gamma_{j}$.
We then observe that this quantity is unstable 
in the sense that it becomes extremely large
for blurry images obtained by convolving
a cartoon image by a mollifier $\vp_\delta=\frac 1{\delta^{2}}\vp(\frac \cdot h)$
as $\delta \to 0$.  In addition, this quantity does not
control the number of subdomains in the partition.
\item
A recent approach proposed in \cite{DPW} defines
anisotropic smoothness through the geometric 
smoothness properties of the {\it level sets} of the function $f$.
In this approach the measure of smoothness 
is not simple to compute directly from $f$ since it involves each of its level sets
and a smoothness measure of their local parametrization.
\end{enumerate}

The results of \cite{BBLS,CSX,M} describe 
the $L^p$-error of piecewise linear interpolation by 
an optimally adapted triangulation of at most $N$ elements,
when $f$ is a $C^2$ function of two variables. This error is defined
as
$$
\sigma_N(f)_p:=\inf_{\#(\cT)\leq N}\|f-I_{\cT}f\|_{L^p}.
$$
It is shown in \cite{BBLS} for $p=\infty$
and in \cite{M} for all $1\leq p\leq \infty$ that
that 
\be
\limsup_{N\to +\infty} N\sigma_N(f)_p\leq CA_p(f),
\label{limsupest}
\ee
where $C$ is an absolute constant and
\be
A_p(f):=\|\sqrt {|{\rm det}(d^2f)|}\|_{L^\tau},\;\; \frac 1\tau:= 1+\frac 1 p.
\label{nonlinnorm}
\ee
Moreover, this estimate is known to be optimal in the sense that
$\liminf_{N\to +\infty} N\sigma_N(f)_p\geq cA_p(f)$ also holds, under some mild restriction
on the class of triangulations in which one selects the optimal one.
These results are extended in \cite{M} to the case of
higher order finite elements and space dimension $d>2$,
for which one can identify similar measures $f\mapsto A(f)$ governing the
convergence estimate. Such quantities thus constitute natural 
candidates to measure anisotropic smoothness properties. Note
that $A_p(f)$ is not a semi-norm due to the presence of the determinant
in \iref{nonlinnorm}, and in particular the
quasi-triangle inequality $A_p(f+g)\leq C(A_p(f)+A_p(g))$ 
does not hold even with $C>1$.

Our paper is organized as follows.
We begin in \S 2 by a brief account of the
available estimates on anisotropic finite element 
interpolation, and we recall in particular the argument that
leads to \iref{limsupest} with the quantity $A_p(f)$ defined by \iref{nonlinnorm}.

Since $A_p(f)$ is not a norm, we cannot associate
a linear smoothness space to it by a standard completion
process. We are thus facing a difficulty in extending
the definition of $A_p(f)$ to functions which are not
$C^2$-smooth and in particular to cartoon images
such as in item 3 above. Since we know from \iref{triN}
that for such cartoon images the $L^2$ error of adaptive
piecewise linear interpolation decays like $N^{-1}$, we would expect that the
quantity 
$$
A_2(f)=\|\sqrt {|{\rm det}(d^2f)|}\|_{L^{2/3}},
$$
corresponding to the case $p=2$ can be properly defined for
piecewise smooth functions. We address this difficulty in 
\S 3 by a regularization process:
if $f$ is a cartoon image we introduce its regularized version
\be
f_\delta:=f*\vp_{\delta},
\ee
where $\vp_\delta=\frac 1{\delta^{2}}\vp(\frac \cdot h)$ 
is a standard mollifier. Our
main result is the following: for any cartoon image $f$ of $C^2-C^2$ type, the quantity $A_2(f_\delta)$
remains uniformly bounded as $\delta\to 0$ and one has
\be
\lim_{\delta \to 0} A_2(f_\delta)^{2/3}= \sum_{i=1}^k \int_{\Omega_i}\left |\sqrt {|{\rm det}(d^2f)|}\right |^{2/3}
+C(\vp)\sum_{j=1}^l \int_{\Gamma_j} |[f](s)|^{2/3} |\kappa(s)|^{1/3} ds,
\label{lima2}
\ee
where $[f](s)$ and $\kappa(s)$ respectively denote
the jump of $f$ and the curvature of $\Gamma_j$ at the point $s$,
and where $C(\vp)$ is a constant that only depends on the choice of the mollifier.
This constant can be shown to be uniformly bounded by below for the class of
radially decreasing mollifiers. 
This result reveals that $A_2(f)$ is stable under regularization of cartoon images
(in contrast to the measure of smoothness described in item 3 above). 
We also discuss the behaviour of $A_p$ when $p\neq 2$.

These results lead us in \S 4 to a comparison
between the quantity $A_2(f)$ and the total variation $\TV(f)$. We also
make some remarks on the existing links between the limit expression
in \iref{lima2} and classical results on adaptive approximation
of curves, as well as
with operators of affine-invariant image processing which also
involve the power $1/3$ of the curvature. 

We devote \S 5 to numerical tests performed on cartoon images that illustrate the validity 
of our results and we describe in \S 6 the extension of our results
to finite elements of higher degree and higher space dimensions.
Concluding remarks and perspectives of our work are given in \S 7.

\section{Anisotropic finite element approximation}

A standard estimate in finite element approximation states that 
if $f\in W^{2,p}(\Omega)$ then
$$
\|f-I_{\cT_h} f\|_{L^p} \leq C  h^{2}\|d^2f\|_{L^{p}},
$$
where $\cT_h$ is a triangulation of mesh size
$h:=\max_{T\in\cT_h} {\rm diam}(T)$. If we restrict our attention
to a family {\it quasi-uniform} triangulations, $h$ is linked with 
the complexity $N:=\#(\cT_h)$ according to
$$
C_1 h^{-2} \leq N \leq  C_2 h^{-2}
$$
Therefore, denoting by $\sigma^{\rm unif}_N(f)_{L^p}$ the $L^p$ approximation
by quasi-uniform triangulations of cardinality $N$, we can re-express the above
estimate as
\be
\label{uniferror}
\sigma^{\rm unif}_N(f)_{L^p} \leq CN^{-1}\|d^2f\|_{L^p}.
\ee
In order to explain how this estimate can be improved
when using adaptive partitions, we first give some
heuristic arguments which are based on the assumption that on each triangle $T$
the relative variation of $d^2f$ is small so that it can be considered
as constant over $T$, which means that 
$f$ coincides with a quadratic function $q_T$ on each $T$.
Denoting by $I_T$ the local interpolation operator
on a triangle $T$ and by $e_{T}(f)_p:=\|f-I_{T}f\|_{L^p(T)}$
the local $L^p$ error, we thus have according to
this heuristics
$$
\|f-I_{\cT}f\|_{L^p} =\(\sum_{T\in \cT} e_T(f)_p^p\)^{\frac 1 p}=\(\sum_{T\in \cT} e_T(q_T)_p^p\)^{\frac 1 p}
$$
We are thus led to study the local interpolation error $e_T(q)_p$ when $q\in \P_2$
is a a quadratic polynomial. Denoting by $\bq$ the homogeneous part of $q$, we
remark that
$$
e_T(q)_p=e_T(\bq)_p.
$$
We optimize the shape of $T$ with respect to the quadratic form $\bq$
by introducing a function $K_p$ defined on the space of quadratic forms by
$$
K_p(\bq):=\inf_{|T|=1}e_T(\bq)_p,
$$
where the infimum is taken among all triangles of area $1$. 
It is easily seen that $e_T(\bq)_p$ is invariant
by translation of $T$ and so is therefore the minimizing triangle
if it exists. By homogeneity, it is also easily seen
that 
$$
\inf_{|T|=a} e_T(\bq)_p=a^{\frac 1 \tau}K_p(\bq),\;\; \frac 1 \tau:=\frac 1 p+1,
$$
and that the minimizing triangle of area $a$ is obtained by rescaling the minimizing
triangle of area $1$ if it exists. Finally, it is easily seen that if $\vp$ is an invertible
linear transform 
$$
K_p(\bq \circ \vp)=|{\rm det} (\vp)|K_p(\bq),
$$
and that the minimizing triangle of area $|{\rm det} (\vp)|^{-1}$ for $\bq \circ \vp$
is obtained by application of $\vp^{-1}$ to the minimizing triangle of 
area $1$ for $\bq$ if it exists. If ${\rm det} (\bq)\neq 0$, there exists 
a $\vp$ such that $\bq \circ \vp$ is either $x^2+y^2$ or $x^2-y^2$ up to 
a sign change, and we have $|{\rm det} (\bq)|=|{\rm det} (\vp)|^{-2}$.
It follows that $K_p(\bq)$ has the simple form
\be
K_p(\bq)=\sigma  |{\rm det} (\bq)|^{1/2},
\label{Kpexp}
\ee
where $\sigma$ is a constant 
equal to $K_p(x^2+y^2)$ 
if ${\rm det} (\bq)>0$ and to  $K_p(x^2-y^2)$ if ${\rm det} (\bq)<0$. 
One easily checks that this equality
also holds when ${\rm det} (\bq)= 0$ in which case $K_p(\bq)=0$.

Assuming that the triangulation $\cT$ is
such that all its triangles $T$ have optimized shape in the above sense
with respect to the quadratic form $\bq_T$ associated with $q_T$, we thus have
for any triangle $T\in\cT$
$$
e_T(f)_p=e_T(\bq_T)_p=|T|^{\frac 1 \tau}K_p(\bq_T)=\left\|K_p\(\frac {d^2f}2\)\right\|_{L^\tau(T)}.
$$
since we have assumed $\frac {d^2f} 2=\bq_T$ on $T$. In order to optimize the trade-off
between the global error and the complexity $N=\#(\cT)$,
we apply the principle of {\it error equidistribution}: the triangles $T$ have
area such that all errors $e_T(\bq_T)_p$ are equal
i.e. $e_T(\bq_T)_p= \eta$ for some $\eta>0$ independent of $T$.
It follows that
$$
N\eta^\tau\leq \left\|K_p\(\frac {d^2f}2\)\right\|_{L^\tau(\Omega)}^\tau,
$$
and therefore
$$
\sigma_N(f)_p \leq \|f-I_\cT f\|_{L^p} \leq N^{1/p} \eta \leq\left\|K_p\(\frac {d^2f}2\)\right\|_{L^\tau(\Omega)}N^{-1},
$$
which according to \iref{Kpexp} implies 
\be
\sigma_N(f)_p \leq CN^{-1}A_p(f),
\label{false}
\ee
with $A_p$ defined as in \iref{nonlinnorm}.

The estimate \iref{false} is too optimistic to be correct:
if $f$ is a univariate function then $A_p(f)=0$ while $\sigma_N(f)_p$
may not vanish. In a rigorous 
derivation such as in \cite{BBLS}, one observes that
if $f\in C^2$, the replacement of $d^2f$ by a constant over $T$
induces an error which becomes negligible only when the
triangles are sufficiently small, and therefore a
correct statement is that for any $\e>0$ there exists 
$N_0=N_0(f,\e)$ such that 
\be
\sigma_N(f)_p \leq N^{-1}\left(\left\|K_p\(\frac {d^2f}2\)\right\|_{L^\tau(\Omega)}+\e\right),
\label{true}
\ee
for all $N\geq N_0$, i.e.
\be
\limsup_{N\to +\infty} N\sigma_N(f)_p\leq \left\|K_p\(\frac {d^2f}2\)\right\|_{L^\tau(\Omega)}
\label{limsupest1}
\ee
which according to \iref{Kpexp} implies \iref{limsupest}.

\section{Piecewise smooth functions and images}

As already observed, the quantities $A_p(f)$ are well defined for
functions $f\in C^2$, but 
we expect that they should in some sense also be well defined for
functions representing $C^2-C^2$ ``cartoon images'' when $p\leq 2$. We first give 
a precise definition of such functions.

\begin{definition} 
\label{defcartoon}
A cartoon function on an open set $\Omega$ is a function almost everywhere of the form 
$$
f=\sum_{1\leq i\leq k} f_i \Chi_{\Omega_i},
$$
where the $\Omega_i$ are disjoint open sets with piecewise $C^2$ boundary, no cusps (i.e. satisfying an interior and exterior cone condition), and such that $\overline \Omega = \cup_{i=1}^k \overline \Omega_i$, and where for each $1\leq i\leq k$, the function $f_i$ is $C^2$ on $\overline \Omega_i$. 
\end{definition}

Let us consider a fixed cartoon function $f$ on an open polygonal domain $\Omega$ (i.e.
$\Omega$ is such that $\overline\Omega$ is a closed polygon)
associated with a decomposition $(\Omega_i)_{1\leq i\leq k}$. We define 
\be
\Gamma :=\bigcup_{1\leq i\leq k} \partial \Omega_i,
\ee
the union of the boundaries of the $\Omega_i$. The above definition
implies that $\Gamma$ is the disjoint union of a finite set of points $\cP$ and a finite number of open curves $(\Gamma_i)_{1\leq i\leq l}$.
$$
\Gamma = \(\bigcup_{1\leq i\leq l} \Gamma_i\) \cup \cP.
$$
Furthermore for all $1\leq i < j \leq l$, we may impose that $\overline \Gamma_i\cap \overline \Gamma_j \subset \cP$ (this may be ensured by a splitting of some of the $\Gamma_i$ if necessary).

We now consider the piecewise linear interpolation
$I_{\cT_N}f$ of $f$ on a triangulation $\cT_N$ of cardinality $N$. 
We distinguish two types of elements of $\cT_N$. A 
triangle $T\in \cT_N$ is called ``regular'' if $T\cap \Gamma=\emptyset$, 
and we denote the set of such triangles by $\cT_N^r$. 
Other triangles are called ``edgy'' and their set is denoted by $\cT_N^e$.
We can thus split $\Omega$ according to
$$
\Omega:=\left(\cup_{T\in \cT_N^r}T\right) \cup \left(\cup_{T\in \cT_N^e}T\right)=\Omega_N^r \cup \Omega_N^e.
$$
We split accordingly the $L^p$ interpolation error into
$$
\|f-I_{\cT_N}f\|_{L^p(\Omega)}^p = \int_{\Omega_N^r} |f-I_{\cT_N}f|^p +\int_{\Omega_N^e} |f-I_{\cT_N}f|^p.
$$
We may use $\cO(N)$ triangles in $\cT_N^e$ and $\cT_N^r$ (for example $N/2$
in each set). Since $f$ has discontinuities along $\Gamma$, the $L^\infty$ interpolation error on 
$\Omega_N^e$ does not tend to zero and $\cT_N^e$ should be chosen so 
that $\Omega_N^e$ has the aspect of a thin layer around $\Gamma$. 
Since $\Gamma$  is a finite union of $C^2$ curves, we can build this layer
of width $\cO(N^{-2})$ and therefore of global area $|\Omega_N^e|\leq C N^{-2}$,
by choosing long and thin triangles in $\cT_N^e$. 
On the other hand, since $f$ is uniformly $C^2$ on $\Omega_N^r$, we may
choose all triangles in $\cT_N^r$ of regular shape 
and diameter $h_T\leq C N^{-1/2}$.
Hence we obtain the following heuristic error estimate, 
for a well designed anisotropic triangulation:
\begin{eqnarray*}
\|f-I_{\cT_N}f\|_{L^p(\Omega)} &=& \(\|f-I_{\cT_N}f\|_{L^p(\Omega_N^r)}^p +\|f-I_{\cT_N}f\|_{L^p(\Omega_N^e)}^p\)^{1/p}\\
&\leq&\( \|f-I_{\cT_N}f\|_{L^\infty(\Omega_N^r)}^p |\Omega_N^r|+  
\|f-I_{\cT_N}f\|_{L^\infty(\Omega_N^e)}^p |\Omega_N^e|\)^{1/p}\\
& \leq & C(N^{-p}+ N^{-2})^{1/p},
\end{eqnarray*}
and therefore
\be
 \|f-I_{\cT_N}f\|_{L^p(\Omega)} \leq C N^{-\min\{1,2/p\}},
 \label{heurist}
 \ee
where the constant $C$ depends on $\|d^2f\|_{L^\infty(\Omega\sm \Gamma)}$, $\|f\|_{L^\infty(\Omega)}$
and on the number, length and maximal curvature of the $C^2$ curves which constitute $\Gamma$.

Observe in particular that the error
is dominated by the edge term $\|f-I_{\cT_N}f\|_{L^p(\Omega_N^e)}$
when $p>2$ and by the smooth term $\|f-I_{\cT_N}f\|_{L^p(\Omega_N^r)}$
when $p<2$. For the critical value $p=2$ the two terms have the same order. 

For $p\leq 2$, we obtain the approximation rate $N^{-1}$ 
which suggests that approximation results such as \iref{limsupest}
should also apply to cartoon functions and
that the quantity $A_p(f)$ should be
finite. We would therefore like to bridge the gap 
between anisotropic approximation of cartoon functions and smooth functions. 
For this purpose, we first need to give a proper meaning
to $A_p(f)$ when $f$ is a cartoon function. This is not 
straightforward, due to the fact that
the product of two distributions has no meaning in general.
Therefore, we cannot define $\det (d^2 f)$ in a distributional sense,
when the coefficients of $d^2 f$ are distributions without sufficient smoothness.
Our approach will rather be based on regularisation. This is additionally justified by the fact that sharp curves of discontinuity are a mathematical idealisation. In real world applications, such as photography, several physical limitations (depth of field, optical blurring) impose a certain level of blur on the edges.

In the following, we consider a fixed
radial nonnegative function $\vp$ of unit integral and supported in the unit ball,
and we define for all $\delta>0$ and $f$ defined on $\Omega$,
\be
\vp_\delta(z) := \frac 1 {\delta^2} \vp\left(\frac z \delta\right) \text{ and } f_\delta = f * \vp_\delta.
\label{defphih}
\ee
Our main result gives a meaning to $A_p(f)$ based on this regularization.
If $f$ is a cartoon function on a set $\Omega$, 
and if $x\in \Gamma\sm\cP$, we denote by $[f](x)$ the jump of $f$
at this point.
We also denote $\bt(x)$ and $\bn(x)$ the unit tangent and normal
vectors to $\Gamma$ at $x$ oriented in such way that $\det(\bt,\bn)=+1$,
and by $\kappa(x)$ the curvature
at $x$ which is defined by the relation
$$
\partial_{\bt(x)}\bt(x)=\kappa(x)\bn(x).
$$
For $p\in [1,\infty]$ and $\tau$ defined by $\frac 1 \tau := 1+\frac 1 p$,
we introduce the two quantities
\begin{eqnarray*}
S_p(f) &:= & \|\sqrt{|\det (d^2 f) |}\|_{L^\tau(\Omega\sm \Gamma)}=A_p(f_{|\Omega\sm \Gamma}),\\
E_p(f) &:= & \|\sqrt{|\kappa|} [f]\|_{L^\tau(\Gamma)},\\
\end{eqnarray*}
which respectively measure the ``smooth part'' and the ``edge part'' of $f$. We
also introduce the constant
\be
C_{p,\vp}:=\|\sqrt {|\Phi\Phi'|}\|_{L^\tau(\R)},\;\; \text{ where } \ \Phi(x) := \int_{y\in\R} \vp(x,y)dy.
\label{defchi}
\ee
Note that $f_\delta$ is only properly defined on the set
$$
\Omega^\delta:=\{z\in \Omega\; ; \; B(z,\delta)\subset \Omega\},
$$
and therefore, we define
$A_p(f_\delta)$ as the $L^\tau$ norm of $\sqrt{|\det(d^2f_\delta)|}$ on this set.

\begin{theorem}
\label{threg}
For all cartoon functions $f$, the quantity $A_p(f_\delta)$ behaves as follows:
\begin{itemize}
\item If $p<2$, then
$$
\lim_{\delta\to 0} A_p(f_\delta) = S_p(f).
$$
\item If $p=2$, then $\tau=\frac 2 3$ and 
\be
\lim_{\delta\to 0} A_2(f_\delta) = \(S_2(f)^\tau+ E_2(f)^\tau C_{2,\vp}^\tau\)^{\frac 1 \tau}.
\label{Acartoon}
\ee
\item If $p>2$, then $A_p(f_\delta) \to \infty$ according to
\be
\lim_{\delta\to 0}  \delta^{\frac 1 2 - \frac 1 p} A_p(f_\delta)= E_p(f) C_{p,\vp}.
\label{equivAphih}
\ee
\end{itemize}
\end{theorem}

\begin{remark}
This theorem reveals that as $\delta\to 0$, the 
contribution of the neighbourhood of $\Gamma$ to $A_p(f_\delta)$
is negligible when $p<2$ and 
dominant when $p>2$, which was already remarked in
the heuristic computation leading to {\rm \iref{heurist}}.
\end{remark}

\begin{remark}
It seems to be possible to eliminate the ``no cusps'' condition
in the definition of cartoon functions, while still retaining the validity of
this theorem. It also seems possible to take the more natural choice 
$\vp(z) = \frac 1 \pi e^{-\|z\|^2}$, which is not compactly supported.
However, both require higher technicality in the proof
which we avoid here.
\end{remark}

Before attacking the proof of Theorem \ref{threg},
we show below that the constant $C_{p,\vp}$ involved in 
the result for $p\geq 2$ is uniformly bounded by below
for a mild class of mollifiers.

\begin{proposition}
Let $\vp$ be a radial and positive function supported on the unit ball
such that $\int \vp=1$ and that $\vp(x)$ decreases as $|x|$ increases. For any $p\geq 2$ we have 
$$
C_{p, \vp} \geq \frac 2 \pi \left(\frac 4 {\tau+2}\right)^{\frac 1 \tau}.
$$
and this lower bound is optimal. 
There is no such bound if $p<2$, but note that Theorem \ref{threg} does not involve $C_{p, \vp}$ for $p<2$.
\end{proposition}

\proof
Let $D$ be the unit disc of $\R^2$. We define a non smooth mollifier $\psi$ and a function $\Psi$ as follows
$$
\psi := \frac {\Chi_D} \pi \ \text{ and } \ \Psi(x) :=  \int_\R\psi(x,y) dy.
$$
One easily obtains that
$
\Psi(x) = \frac 2 \pi \sqrt{1-x^2} \Chi_{[-1,1]}(x)\text{ and } \Psi'(x) = \frac {-2x}{\pi \sqrt{1-x^2}} \Chi_{[-1,1]}(x), 
$
hence 
$$
\Psi(x) \Psi' (x)=  \frac {-4x} {\pi^2} \Chi_{[-1,1]}(x).
$$
For all $\delta>0$ we define $\psi_\delta := \delta^{-2} \psi(\delta^{-1} \cdot)$, and  $\Psi_\delta(x) :=  \int_\R\psi_\delta(x,y) dy$. Similarly we obtain 
$$
\Psi_\delta(x) \Psi'_\delta (x) = \frac {-4x} {\pi^2}  \Chi_{[-\delta,\delta]} \delta^{-4}.
$$
Hence 
$$
C_{p,\psi_\delta} = \left\|\sqrt{\Psi_\delta(x) |\Psi'_\delta (x)|}\right\|_{L^\tau(\R)} = \frac 2 {\pi \delta^2} \left(\int_{-\delta}^\delta |x|^{\frac\tau 2} dx\right)^{\frac 1 \tau} 
= \frac 2 {\pi \delta^2} \left(\frac{2 \delta^{\frac \tau 2 +1}}{\frac \tau 2 +1}\right)^{\frac 1 \tau}
= \frac 2 \pi \left(\frac 4 {\tau+2}\right)^{\frac 1 \tau} \delta^{\frac 1 p -\frac 1 2}.
$$
Note that
\be
\label{ineqPsiDelta}
\text{ If } p\geq 2 \text{ and }\delta \in (0,1] \text{ then }  C_{p,\psi_\delta} \geq C_{p,\psi} =   \frac 2 \pi \left(\frac 4 {\tau+2}\right)^{\frac 1 \tau}.
\ee
The mollifier $\vp$ of interest is radially decreasing, has unit integral and is supported on the unit ball. It follows that there exists a Lebesgue measure $\mu$ on $(0,1]$, of mass $1$, such that 
$$
\vp = \int_0^1\psi_\delta \ d\mu(\delta).
$$
Hence
$
\Phi(x) := \int_\R \vp(x, y) dy = \int_0^1\Psi_\delta(x) \ d\mu(\delta)
$,
 for any $x\in \R$.
Since $s\mapsto s^\tau$ is concave on $\R_+$ when $0<\tau \leq 1$, we obtain
$$
\Phi(x)^\tau = \left(\int_0^1 \Psi_\delta(x) \ d\mu(\delta)\right)^\tau \geq \int_0^1 \Psi_\delta(x)^\tau \ d\mu(\delta)
$$
Similarly, since the sign of $\Psi'_\delta(x)$ is independent of $\delta$, 
$
|\Phi'(x)|^\tau = \left(\int_0^1 |\Psi_\delta'(x)| \ d\mu(\delta)\right)^\tau \geq \int_0^1 |\Psi_\delta'(x)|^\tau \ d\mu(\delta).
$
Applying the Cauchy-Schwartz inequality we obtain
$$
\sqrt{\Phi(x) |\Phi'(x)|}^{\, \tau} \geq 
\sqrt{\left(\int_0^1 \Psi_\delta(x)^\tau \ d\mu(\delta)\right)\left(\int_0^1 |\Psi_\delta'(x)|^\tau \ d\mu(\delta)\right)} 
\geq \int_0^1 \sqrt{\Psi_\delta(x) |\Psi_\delta'(x)|}^{\, \tau} \ d\mu(\delta)
$$
Eventually we conclude using the previous equation and \iref{ineqPsiDelta} as follows
$$
C_{p,\vp}^\tau = \int_\R \sqrt{\Phi |\Phi'|}^{\,\tau}\\
 \geq \int_0^1 \left(\int_\R \sqrt{\Psi_\delta |\Psi_\delta'|}^{\,\tau}\right) d\mu(\delta)\\
=  \int_0^1  C_{p,\psi_\delta}^\tau d\mu(\delta) \geq  C_{p,\psi}^\tau,
$$
which concludes the proof of this lemma.
\sq

The rest of this section is devoted to the proof of Theorem \ref{threg}. Since
it is rather involved, we split its presentation into several main steps.
\nl
\nl
{\bf Step 1: decomposition of $A_p(f_\delta)$.}
Using the notation $K(M):= \sqrt{|\det M|}$, we can write
\be
\label{tointegrate2d}
A_p(f_\delta)^\tau = \int_{\Omega^\delta} K(d^2 f_\delta)^\tau.
\ee
We decompose this quantity based on a partition of
$\Omega^\delta$ into three subsets 
$$
\Omega^\delta=\Omega_\delta \cup \Gamma_\delta \cup \cP_\delta.
$$
The first set $\Omega_\delta$ corresponds to the
{\it smooth part}:
$$
\Omega_\delta :=  \bigcup_{1\leq i\leq k} \Omega_{i,\ssdelta}, \;\;{\rm where}\;\; 
\Omega_{i,\ssdelta} := \{z\in \Omega_i\sep d(z,\Omega\sm\Omega_i)>\delta\}.
$$
Note that $\Omega_\delta$ is strictly contained in $\Omega^\delta$.
The second set corresponds to the {\it edge part}:  we first define
$$
 \Gamma_\delta^0 := \bigcup_{1\leq j\leq l} \Gamma_{j,\ssdelta}^0, \;\;{\rm where}\;\; 
\Gamma_{j,\ssdelta}^0 := \{z\in \Gamma_j\sep d(z,\Gamma\sm\Gamma_j)>2\delta\},
$$
and then set
$$
\Gamma_\delta := \bigcup_{1\leq j\leq l} \Gamma_{j,\ssdelta}\;\;{\rm where}\;\; 
\Gamma_{j,\ssdelta} :=  \{z\in \Omega \sep d(z,\Gamma)< \delta \text{ and } \pi_\Gamma(z)\in \Gamma_{j,\ssdelta}^0\} 
$$
where $\pi_\Gamma(z)$ denotes the point of $\Gamma$ which is the closest to $z$.
The third set corresponds the {\it corner part}:
$$
\cP_\delta :=  \Omega^\delta \sm (\Omega_\delta\cup\Gamma_\delta).
$$
The measures of the sets $\Gamma_\delta$ and $\cP_\delta$ tends to
$0$ as $\delta\to 0$, while $|\Omega_\delta|$ tends to $|\Omega|$. More
precisely, we have
$$
 |\Gamma_\delta|\leq C\delta\;\;{\rm and}\;\; |\cP_\delta|\leq C  \delta^2
$$
where the last estimate exploits the ``no cusps'' property of the cartoon function.
We analyze separately the contributions of these three sets to 
\iref{tointegrate2d}.
\nl
\nl
{\bf Step 2: Contribution of the smooth part $\Omega_\delta$.}
The contribution of $\Omega_\delta$ to the integral \iref{tointegrate2d} is easily measured. Indeed, let us define
$$
Q_\delta(z):= \left\{
\begin{array}{cc}
K(d^2 f_\delta (z))^\tau & \text{ if } z\in \Omega_\delta,\\
0 & \text{ otherwise.}
\end{array}
\right.
$$
then we have pointwise convergence $Q_\delta(z)\to K(d^2 f(z))^\tau$ on $\Omega\sm \Gamma$.
Since the  $\delta$-neighbourhood of $\Omega_\delta$ is included in $\Omega\sm\Gamma$, we have 
$$
\|d^2 (f*\vp_\delta) \|_{L^\infty(\Omega_\delta)} = \|(d^2 f)*\vp_\delta \|_{L^\infty(\Omega_\delta)}\leq \|d^2 f\|_{L^\infty(\Omega\sm \Gamma)} \|\vp_\delta\|_{L^1} = \|d^2 f\|_{L^\infty(\Omega\sm \Gamma)} 
$$
Since $K$ is $1$-homogeneous and continuous, we have
$$
K(d^2 f_\delta)\leq C_K\|d^2 f\|_{L^\infty(\Omega\sm \Gamma)}, \;\; C_K:=\max_{\|M\|=1}K(M),
$$
and we conclude by dominated convergence that 
$$
\lim_{\delta \to 0}\int_{\Omega_\delta} K(d^2 f_\delta)^\tau = \lim_{\delta\to 0} \int_{\Omega \sm \Gamma} Q_\delta =  \int_{\Omega \sm \Gamma} K(d^2 f)^\tau.
$$
\nl
{\bf Step 3: Contribution of the corner part $\cP_\delta$.}
We only need a rough upper estimate of the contribution of $\cP_\delta$ 
to the integral \iref{tointegrate2d}.
We observe that 
$$
\|d^2 (f*\vp_\delta) \|_{L^\infty(\Omega)} = \|f*(d^2  \vp_\delta)\|_{L^\infty(\Omega)}\leq \|f\|_{L^\infty(\Omega)} \|d^2 \vp_\delta\|_{L^1(\R^2)} =   \frac M {\delta^2},
$$
where $M:=\|f\|_{L^\infty(\Omega)} \|d^2 \vp\|_{L^1(\R^2)}$.
It follows that 
$$
\int_{\cP_\delta} K(d^2 f_\delta)^\tau \leq |\cP_\delta| \left(C_K \left(\frac M {\delta^2}\right)\right)^\tau \leq C  \delta^{2-2 \tau}.
$$
If $\tau<1$, this quantity tends to $0$ and is therefore negligible compared to the 
contribution of the smooth part. If $\tau=1$, which corresponds to $p=\infty$, our
further analysis shows that the contribution of the edge part tends
to $+\infty$, and therefore the contribution of the corner part is always negligible.
\nl
\nl
{\bf Step 4: Contribution of the edge part $\Gamma_\delta$.} {
This step is the main difficulty of the proof. We make a key use
of an asymptotic analysis of $f_\delta$ on $\Gamma_\delta$, which relates
its second derivatives to the jump $[f]$ and the curvature $\kappa$ as $\delta\to 0$. 
We first define for all $\delta>0$ the map
$$
\begin{array}{ccc}
U_\delta : \Gamma\sm\cP \times [-1,1] &\to& \Omega\\
(x,u) &\mapsto & x+\delta u \bn (x).
\end{array}
$$
We notice that according to our definitions,
for $\delta$ small enough, the map $U_\delta$ induces a diffeomorphism 
between $\Gamma_\delta^0\times [-1,1]$ and $\Gamma_\delta$, such
that $\pi_\Gamma(U_\delta(x,u))=x$ and $d(U_\delta(x,u),\Gamma)=|U_\delta(x,u)-x|=\delta |u|$.
We establish asymptotic estimates on the
second derivatives of $f_\delta$ which have the following form:
\begin{eqnarray}
\label{estimnn2d}
\left |\partial_{\n,\n} f_\delta(z) -  \frac 1 { \delta^2} [f](x) \Phi'(u)\right| &\leq & \frac C { \delta}\\ 
\label{estimnt2d}
|\partial_{\n,\t} f_\delta(z)| &\leq &  \frac C \delta\\
\label{estimtt2d}
\left |\partial_{\t,\t} f_\delta(z) +  \frac 1 \delta [f](x) \kappa(x)\Phi(u)\right| &\leq &\frac {\omega(\delta)} \delta
\end{eqnarray}
where $\lim_{\delta\to 0} \omega(\delta) = 0$ and with the notation
$z=U_\delta(x,u)$. The constant $C$ and the function $\omega$ depend only on $f$.
The proof of these estimates is given in the appendix. As an immediate consequence,
we obtain an asymptotic estimate of $K(d^2 f_\delta)=\sqrt{|\det(d^2 f_\delta)|}$ of the form
\be
\label{estimK2d}
\left|  \delta^{\frac 3 2} K(d^2 f_\delta(z)) - \sqrt{|\kappa(x)|} \, |[f](x)| \sqrt{|\Phi(u) \Phi'(u)|} \right|\leq \omega(\delta),
\ee
where $\lim_{\delta\to 0} \omega(\delta) = 0$, and the function $\omega$ depends only on $f$. 
Using the notations
$$
g_\delta(z):=\delta^{\frac 3 2} K(d^2 f_\delta(z)),\;\;  \lambda(x):=\sqrt{|\kappa(x)|}\, |[f](x)|,\;\; \mu(u):=\sqrt{|\Phi(u) \Phi'(u)|},
$$
we thus have
\be
|g_\delta(z) - \lambda(x) \mu(u)|\leq \omega(\delta), 
\label{glambmu}
\ee
for all $x\in \Gamma^0_\delta$, $u\in [-1,1]$ and $\delta>0$
sufficiently small, with $z=U_\delta(x,u)$. We claim that
for any continuous functions $(g_\delta,\lambda,\mu)$ satisfying \iref{glambmu},
we have for any $\tau>0$, 
\be
\lim_{\delta\to 0}  \delta^{-1} \int_{\Gamma_\delta} g_\delta^\tau =\int_{\Gamma} \lambda(x)^\tau dx  \int_{-1}^1 \mu(u)^\tau du,
\label{claim}
\ee
which is in our case equivalent to the estimate
\be
\lim_{\delta\to 0}  \delta^{ \frac 3 2 \tau-1} \int_{\Gamma_\delta} K(d^2 f_\delta)^\tau = 
 \int_\Gamma |[f]|^\tau |\kappa|^{\tau/2} \int_\R |\Phi \Phi'|^{\tau/2}.
\label{estimgammadelta}
\ee
In order to prove \iref{claim}, we may assume without loss of generality that $\tau =1$
up to replacing $(g_\delta,\lambda,\mu)$  by $(g_\delta^\tau,\lambda^\tau,\mu^\tau)$.
We first express the jacobian matrix of $U_\delta$
using the bases $B_1 = ((\bt(x),0),(0,1))$ and $B_2 =(\bt(x),\bn(x))$ 
for the tangent spaces of $\Gamma\times [-1,1]$ and $\Omega$.
This gives the expression
$$
[d U_\delta(x,u)]_{B_1,B_2} = \left(
\begin{array}{c|c}
1- \delta u\kappa(x)& 0\\
\hline
0 & \delta
\end{array}
\right)
$$
and therefore $|\det(dU_\delta(x,u))| = \delta -\delta^2 u\kappa(x)$. Since
$B_1$ and $B_2$ are orthonormal bases, this quantity is the jabobian of $U_\delta$ at $(x,u)$,
and therefore
$$
\int_{\Gamma_\delta} g_\delta =\delta\int_{\Gamma_\delta^0 \times [-1,1]} g_\delta(x+\delta u \bn (x)) 
(1 -\delta u\kappa(x)) dx\, du
$$
Combining with \iref{glambmu}, and using dominated
convergence we obtain \iref{claim}.
\nl
\nl
{\bf Step 5: summation of the different contributions.}
Summing up the contributions of $\Omega_\delta$, $\cP_\delta$ and $\Gamma_\delta$,
we reach the estimate
\begin{eqnarray*}
\int_\Omega K(d^2 f_\delta)^\tau &=&  \int_{\Omega_\delta} K(d^2 f_\delta)^\tau+ \int_{\cP_\delta} K(d^2 f_\delta)^\tau+ \int_{\Gamma_\delta} K(d^2 f_\delta)^\tau\\
&=& \left(\int_{\Omega\sm\Gamma} K(d^2 f)^\tau +\e_1(\delta) 
\right) + B(\delta) \delta^{2-2\tau}+  \delta^{1-\frac 3 2 \tau} \left( \int_\Gamma |[f]|^\tau |\kappa|^{\tau/2}\int_\R |\Phi \Phi'|^{\tau/2} +\e_2(\delta)\right),\\
&=& (S_p(f)^\tau +\e_1(\delta)) + B(\delta) \delta^{2-2\tau}+  \delta^{1-\frac 3 2 \tau} (E_p(f)^\tau C_{p,\phi}^\tau+\e_2(\delta)),\\
\end{eqnarray*}
where $\lim_{\delta\to 0}\e_1(\delta)=\lim_{\delta\to 0}\e_2(\delta)=0$ and $B(\delta)$ is uniformly bounded.
This concludes the proof of Theorem \ref{threg}.

\section{Relation with other works}

Theorem \ref{threg} allows us to extend
the definition of $A_2(f)$ when $f$ is a cartoon
function, according to
\be
A_2(f):= \(S_2(f)^{2/3}+ E_2(f)^{2/3} C_{2,\vp}^{2/3} \)^{3/2}.
\label{addita2}
\ee
We first compare this additive form
with the total variation $\TV(f)$. If $f$ is a
cartoon function,
its total variation has the additive form
\be
\TV(f):=\|\nabla f\|_{L^1(\Omega\sm\Gamma)}+ \|[f]\|_{L^1(\Gamma)},
\label{additv}
\ee
Both \iref{addita2} and \iref{additv}
include a ``smooth term'' and an ``edge term''.
It is interesting to compare the edge term
of $A_2(f)$, which is given by
$$
E_2(f)=\|\sqrt {|\kappa|}[f]\|_{L^{2/3}(\Gamma)},
$$
up to the multiplicative constant $C_{2,\vp}$, with the
one of $\TV(f)$ which is simply the integral of the jump
$$
J(f):=\|[f]\|_{L^{1}(\Gamma)},
$$
Both terms are $1$-homogeneous with the value of the jump of the
function $f$. In particular, if the value of this jump is $1$
(for example when $f$ is the characteristic function of a set of boundary $\Gamma$), we have
\be
E_2(f)=\(\int_\Gamma |\kappa|^{1/3}\)^{3/2},
\label{powercurv}
\ee
while $J(f)$ coincides with the length of $\Gamma$. 
In summary, $A_2(f)$ takes into account the {\it smoothness} of
edges, through their curvature $\kappa$, while 
$\TV(f)$ only takes into account their {\it length}.

Let us now investigate more closely the measure of
smoothness of edges which is incorporated in $A_2(f)$. According to 
\iref{powercurv}, this smoothness is meant in the sense
that the arc length parametrizations of the curves that
constitute $\Gamma$ admit  
second order derivatives in $L^{\frac 1 3}$. In the following,
we show that this particular measure of smoothness is 
naturally related to some known results in two different areas:
adaptive approximation of curves and affine-invariant image
processing.

We first revisit the derivation of the heuristic estimate
\iref{heurist} for the error between a cartoon function
and its linear interpolation on an optimally adapted 
triangulation. In this computation, the contribution of the
``edgy triangles'' was estimated by the area of the layer $\Omega_N^e$ according to
$$
\|f-I_{\cT_N} f\|_{L^p(\Omega^{e}_N)} \leq \|f\|_{L^\infty} |\Omega_N^e|^{1/p}.
$$
Then we invoke the fact that $\Gamma$ is a finite union of $C^2$ curves $\Gamma_j$
in order to build a layer of global area $|\Omega_N^e| \leq CN^{-2}$, 
which results in the case $p=2$ into a contribution to the $L^p$ error
of the order $\cO(N^{-1})$. The area of the layer $\Omega_N^e$ 
is indeed of the same order as the area between the edge $\Gamma$
and its approximation by a polygonal line with $\cO(N)$ segments.

Each of the curves $\Gamma_j$
can be identified to the graph of a $C^2$ function
in a suitable orthogonal coordinate system.
If $\gamma$ is one of these functions, the area
between $\Gamma$
and its polygonal approximation can thus be 
locally measured by the $L^1$ error
between the one-dimensional function
$\gamma$ and a piecewise linear approximation
of this function.
Since $\gamma$ is $C^2$, it is obvious
that it can be approximated by a piecewise linear function on $\cO(N)$ intervals
with accuracy $\cO(N^{-2})$ in the $L^\infty$ norm and therefore
in the $L^1$ norm. However, we may ask whether such a rate
could be achieved under weaker conditions on the smoothness of $\gamma$.
The answer to this question is a chapter of {\it nonlinear} approximation theory
which identifies the exact conditions for a function $\gamma$ to be approximated
at a certain rate by piecewise polynomial functions on adaptive one-dimensional
partitions. We refer to \cite{De} for a detailed
treatment and only state the result which is of interest to us. We say that a function
$\gamma$ defined on a bounded interval $I$ belongs to the approximation
space ${\cal A}^s(L^p)$ if and only if there exists a sequence $(p_N)_{N>0}$ of functions where
each $p_N$ is piecewise affine on a partition of $I$ by $N$ intervals
such that 
$$
\|\gamma-p_N\|_{L^p} \leq CN^{-s}.
$$
For $0<s\leq 2$, it is known that 
$\gamma\in {\cal A}^s(L^p)$ provided that $\gamma\in B^{s}_{\tau,\tau}(I)$ with $\frac 1 \tau:=\frac 1 p+s$,
where $B^s_{\tau,\tau}(I)$ is the standard Besov space
that roughly describes those functions having $s$ derivatives in $L^{\tau}$.
In the case $s=2$ and $p=1$ which is of interest to us, we find $\tau=\frac 1 3$ and therefore $\gamma$ should
belong to the Besov space $B^2_{\frac 1 3,\frac 1 3}(I)$. Note that in our definition 
of cartoon functions, we assume much more than $B^2_{\frac 1 3,\frac 1 3}$ smoothness on $\gamma$,
and it is not clear to us if Theorem \iref{threg} can be derived under this minimal smoothness 
assumption. However it is striking to see that the quantity $E_2(f)$
that is revealed by Theorem \iref{threg} precisely measures the second derivative of the
arc-length parametrization of $\Gamma$ in the $L^{\frac 1 3}$ norm, up to
the multiplicative weight $|[f]|^{2/3}$. Let us also mention
that Besov spaces have been used in \cite{DPW} in order to describe
the smoothness of functions through the regularity of their level sets.
Note that edges and level sets are two distinct concepts, which
coincide in the case of piecewise constant cartoon functions.

The quantity $|\kappa|^{1/3}$ is also encountered in
mathematical image processing, for the design of
simple smoothing semi-groups that respect {\it affine invariance}
with respect to the image. Since these semi-groups should also
have the property of {\it contrast invariance}, they can be defined
through curve evolution operators acting on the level sets
of the image. The simplest curve evolution operator
that respects affine invariance is given by the equation
$$
\frac {d \Gamma}{d t}=-|\kappa|^{1/3}{\rm n},
$$
where ${\rm n}$ is the outer normal, see e.g. \cite{Ca}.
Here the value $1/3$ of the exponent plays a critical
role. The fact that we also find it in $E_2(f)$ suggests
that some affine invariance property also holds for this quantity
as well as for $A_2(f)$. 
We first notice that if $f$ is a compactly supported
$C^2$ function of two variables and $T$ is a 
bijective affine transformation, then with
$\ti f$ such that
$$
f=\ti f\circ T,
$$
we have
the property
$$
d^2 f(z)=L^\trans d^2 \ti f(Tz) L,
$$
where $L$ is the linear part of $T$ and $L^\trans$ its transpose, so that
$$
\sqrt {|\det (d^2f(z))|}=|\det L|\sqrt{|\det (d^2 \ti f(Tz))|}.
$$
By change of variable, we thus find that
\be
A_p(\ti f)=|\det L|^{1/\tau-1}A_p(f)=|\det L|^{1/p}A_p(f).
\label{invap}
\ee
A similar invariance property can be  derived on the interpolation error
$\sigma_N(f)_p=\|f-I_{\cT_N}f\|_{L^p}$ where $\cT_N$ is a triangulation which is optimally adapted
to $f$  in the
sense of minimizing the linear interpolation error in the $L^p$ norm
among all triangulations of cardinality $N$. We indeed remark that
an optimal triangulation for $\ti f$ is then given by applying $T$ to
all elements of $\cT_N$. For such a triangulation 
$\ti \cT_N:=T(\cT_N)$, one has the commutation formula
$$
I_{\cT_N}f=(I_{\ti \cT_N} \ti f)\circ T,
$$
and therefore we obtain by a change of variable that
\be
\sigma_N(\ti f)_p=\|\ti f-I_{\ti\cT_N}\ti f\|_{L^p}=|\det L|^{1/p}\|f-I_{\cT_N}f\|_{L^p}=|\det L|^{1/p}\sigma_N(f)_p.
\label{invsigmap}
\ee
Let us finally show that if $f$ is a cartoon function, then $E_2(f)$ satisfies a similar invariance property
corresponding to $p=2$, namely
\be
E_2(\ti f)=|\det L|^{1/2} E_2(f).
\label{inve2}
\ee
Note that this cannot be derived by arguing that $A_2(f)$ satisfies this
invariance property when $f$ and $\ti f$ are smooth, since we lose
the affine invariance property as we introduce the convolution by $\vp_\delta$: 
we do not have
$$
(f\circ T) * \vp_\delta = (f* \vp_\delta)\circ T.
$$
unless $T$ is a rotation or a translation. 

Let $\Gamma_j$ be one of the $C^2$ pieces of $\Gamma$
and  $\gamma_j: [0,B_i]\to \Omega$ a regular parametrisation of 
$\Gamma_j$. The curvature 
of $\Gamma$ on $\Gamma_j$ at the point $\gamma_j(t)$ is therefore given by
\be
\label{exprKappa}
\kappa(\gamma_j(t))= \frac{\det(\gamma_j'(t),\gamma_j'' (t))}{\|\gamma'_j(t)\|^{3}}
\ee
Since $f=\ti f\circ T$, the discontinuity curves of $\ti f$ are the images of those of $f$ by $T$:
$$
\ti \Gamma_j=T(\Gamma_j).
$$
The curvature of $\ti \Gamma_j$ at the point 
$T(\gamma_j(t))$ is therefore given by 
$$
\ti \kappa(T(\gamma_j(t))) = \frac{\det(L\gamma'_j(t), \ L\gamma''_j(t))}{\|L\gamma'_j(t)\|^3} = \det(L) \frac{\det (\gamma'_j(t),\gamma''_j(t))}{\|L\gamma'_j(t)\|^3} .
$$
This leads us to the relation : 
\be
\label{relKappa}
|\det(L)|^{1/3} |\kappa(\gamma_j(t))|^{1/3} \ \|\gamma'_j(t)\| =  |\ti\kappa(T(\gamma_j(t)))|^{1/3}
\|L\gamma'_j(t)\| ,
\ee
and therefore
\begin{eqnarray*}
\int_{\ti \Gamma_j} |[\ti f]|^{2/3} |\ti \kappa|^{1/3} &=& \int_0^{B_j}  |[\ti f](T(\gamma_j(t)))|^{2/3}
|\ti\kappa(T(\gamma_j(t)))|^{1/3} \ \|L\gamma'_j(t)\| dt,\\
&=& |\det L|^{1/3} \int_0^{B_j}  |[f](\gamma_j(t))|^{2/3} |\kappa(\gamma_j(t))|^{1/3} \  \|\gamma'_j(t)\|dt,\\
&=& |\det L|^{1/3} \int_{ \Gamma_i} |[ f]|^{2/3} |\kappa|^{1/3}.
\end{eqnarray*}
Summing over all $j=1,\cdots,l$ and elevating to the $3/2$ power
we obtain \iref{inve2}.

\section{Numerical tests}

In this section, we validate our previous results
by  numerical tests applied to a simple cartoon image:
the Logan-Shepp phantom. We use a $256\times 256$ pixel
version of this image, with a slight modification which is motivated
further. This image is iteratively smoothed
by the numerical scheme 
\be
\label{numSmoothing}
u_{i,j}^{n+1} = \frac {u_{i,j}^n} 2 + \frac {u_{i+1,j}^n+u_{i-1,j}^n+u_{i,j+1}^n+u_{i,j-1}^n} 8
\ee
This scheme is an explicit discretization of the heat equation.
Formally, as $n$ grows,
$u^n$ is a discretization of 
\be
\label{Gaussian}
u * \vp_{\lambda \sqrt n} \text{ with } \vp_\delta (z) := \frac 1 { \pi\delta^2 } e^{-\frac {\|z\|^2}{\delta^2}},
\ee 
where $u$ stands for the continuous image. The determinant of the hessian 
is discretised by the following $9$-points formula
$$
d_{i,j}^n := (u_{i,j+1}^n-2u_{i,j}^n+u_{i,j-1}^n) (u_{i+1,j}^n-2u_{i,j}^n+u_{i-1,j}^n) - \frac {(u_{i+1,j+1}^n+u_{i-1,j-1}^n-u_{i+1,j-1}^n-u_{i-1,j+1}^n)^2}{16}
$$
For each value of $n$, we then compute the $\ell^\tau$ norm of the array
$\left(\sqrt{|d^n_{i,j}|}\right)$ for $\tau \in [\frac 1 2,1]$, which corresponds to
$p\in [1,\infty]$ with $\frac 1 \tau := 1 + \frac 1 p$. This norm is thus a discretization of the quantity
$$
\left \|\sqrt{|\det(d^2(u*\vp_{\lambda \sqrt n}))|}\right \|_{L^\tau}.
$$
For each value of $n$ we obtain a function $\tau\in  [\frac 1 2,1]\to D_n(\tau) \in\R_+$.

As $n$ grows, three consecutive but potentially overlapping phases appear 
in the behaviour of the functions $D_n$, which are illustrated on Figure 2.

\begin{enumerate}

\item For small $n$, the $9$-points discretisation is not a good approximation of the determinant of the hessian
due to the fact that the pixel discretization is too coarse compared to the 
smoothing width. During this phase, the functions $D_n$ decay rapidly for all values of $\tau$. 

\item For some range of $n$, the edges have been smoothed by the action of \iref{numSmoothing}, 
but the parameter $\lambda \sqrt n$ in \iref{Gaussian} remains rather small. Our previous analysis applies and we observe that $D_n(2/3)$ is (approximately) constant while $D_n(\tau)$ increases for 
$\tau < 2/3$ and decreases for $\tau > 2/3$. 

\item For large $n$, the details of the picture fade and begin to disappear. The picture begins to resemble a constant picture. Therefore the functions $D_n$ decay for all values of $\tau$, and eventually tend to $0$. 

\end{enumerate}

\begin{figure}
	\centering
		\includegraphics[width=5cm,height=5cm]{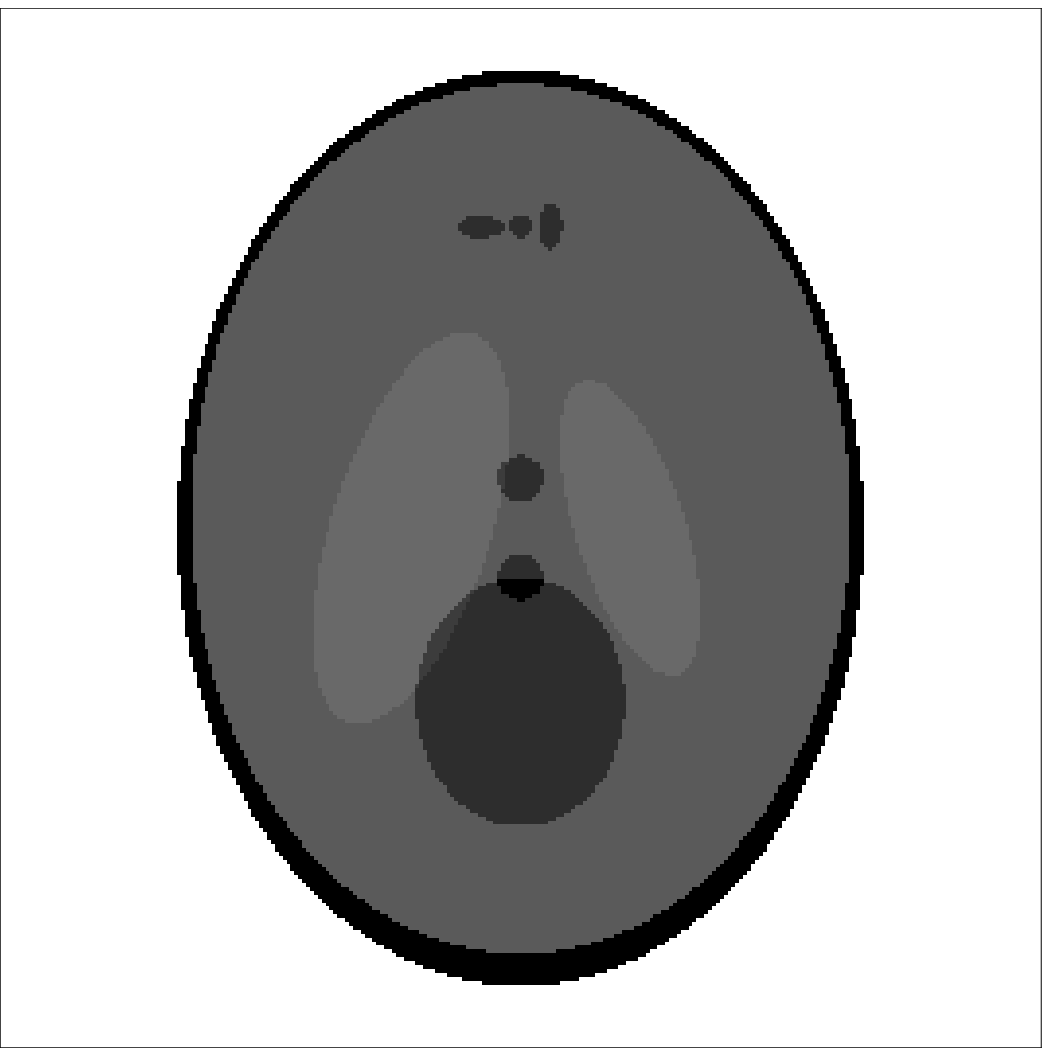}
		\hspace{0.2cm}
		\includegraphics[width=5cm,height=5cm]{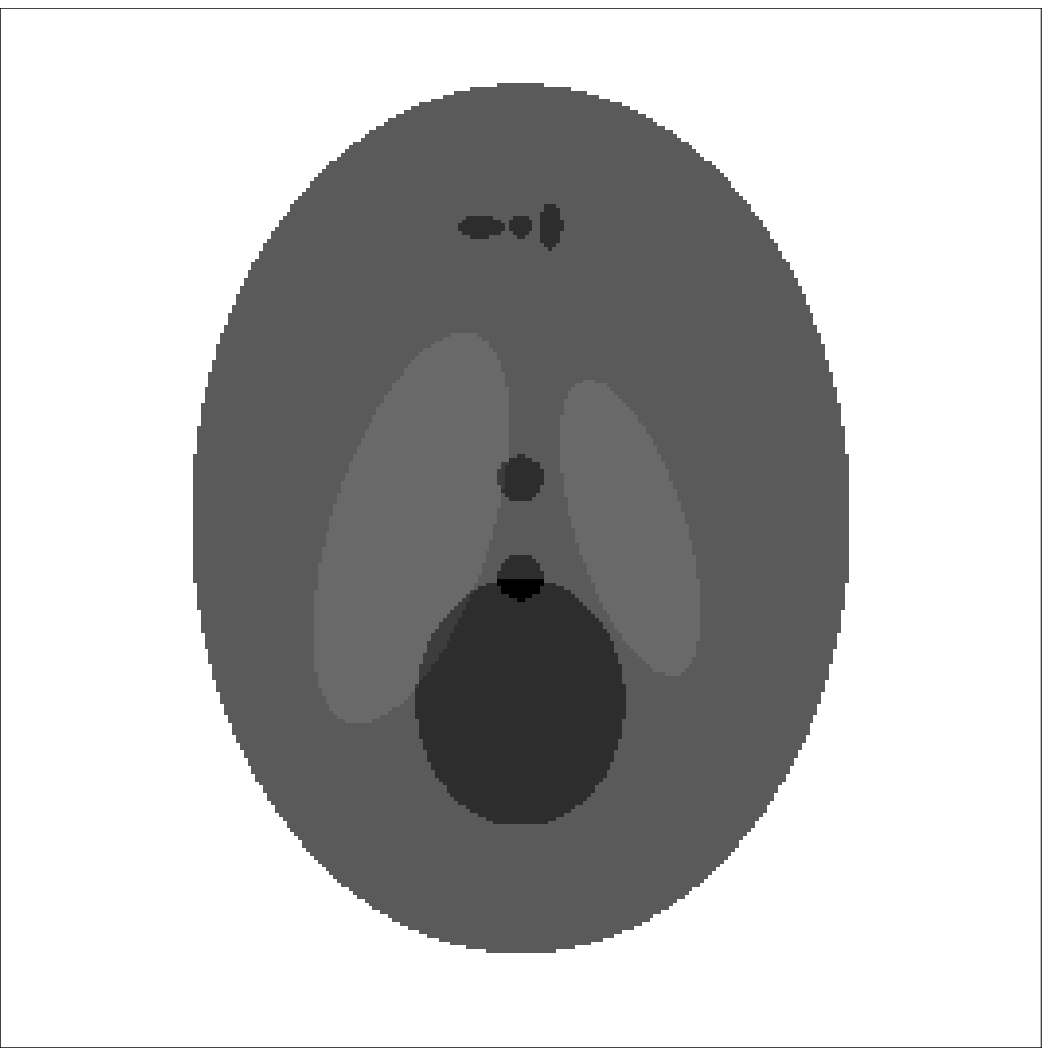}
		\hspace{0.2cm}
		\includegraphics[width=5cm,height=5cm]{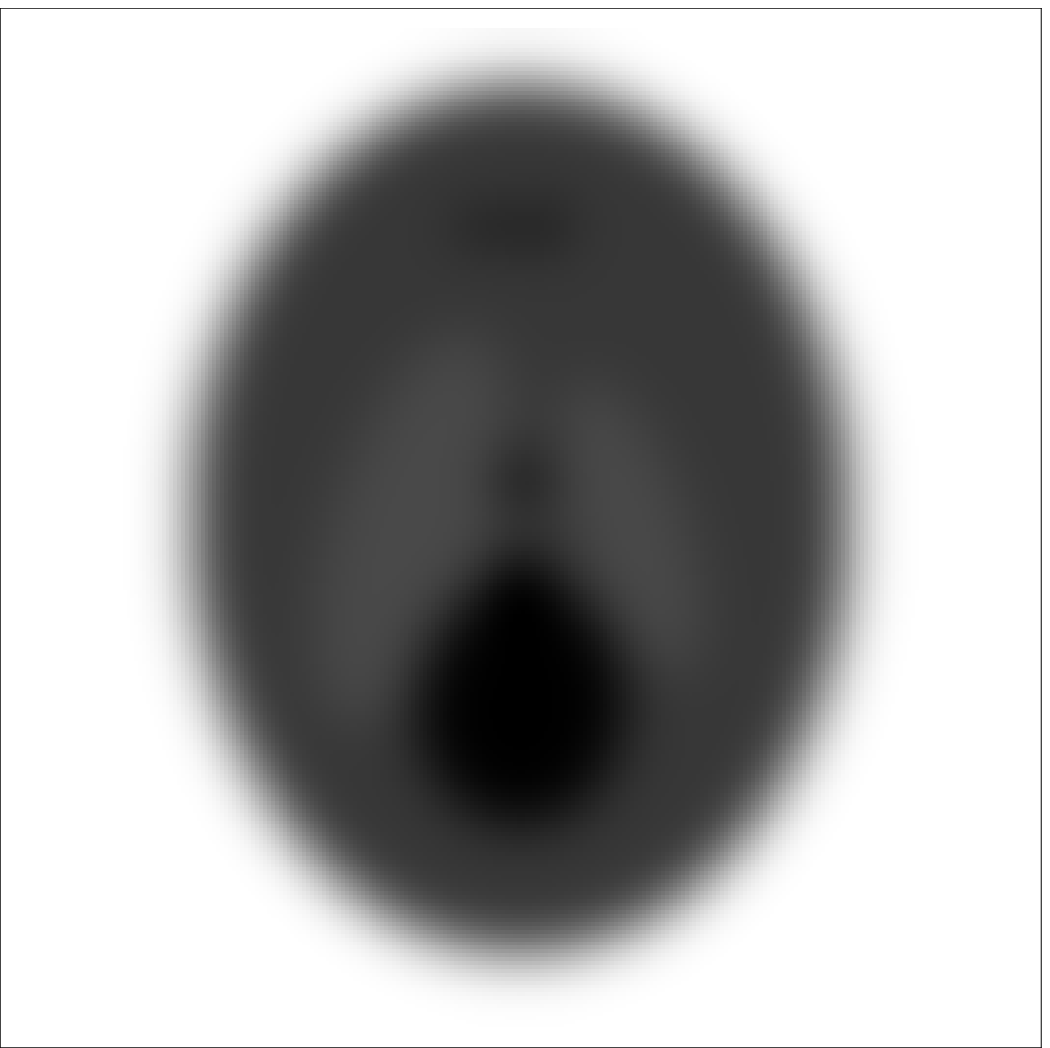}
	\caption{The Logan-Shepp image (left), modified image (center), smoothed image (right).}
\end{figure}

\begin{figure}
	\centering
		\includegraphics[width=5cm,height=4cm]{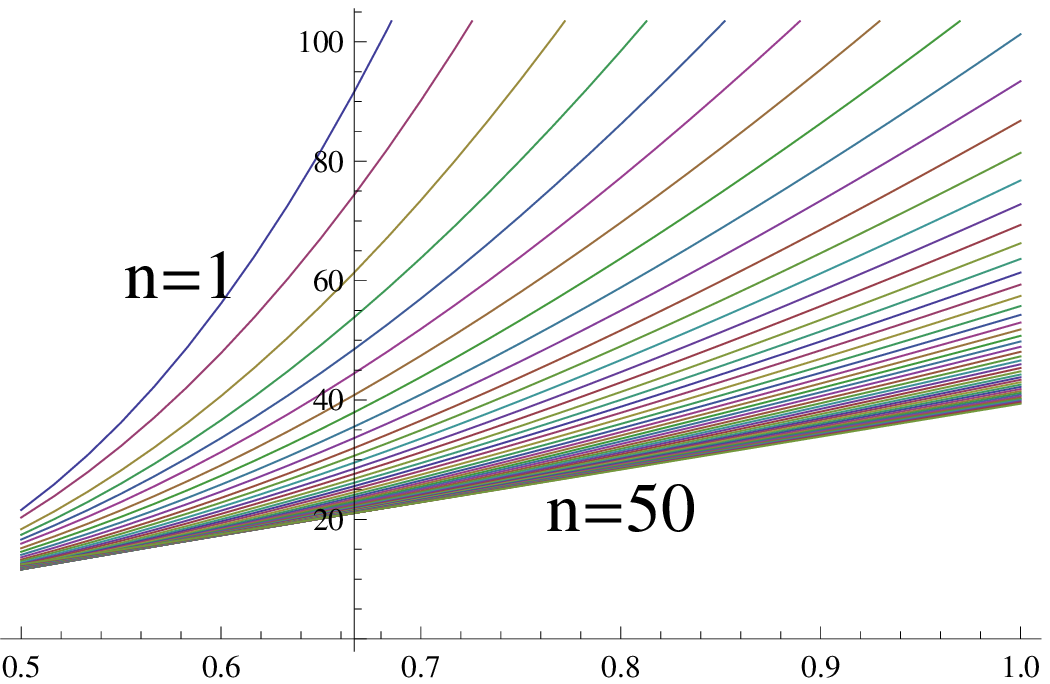}
		\hspace{0.2cm}
		\includegraphics[width=5cm,height=4cm]{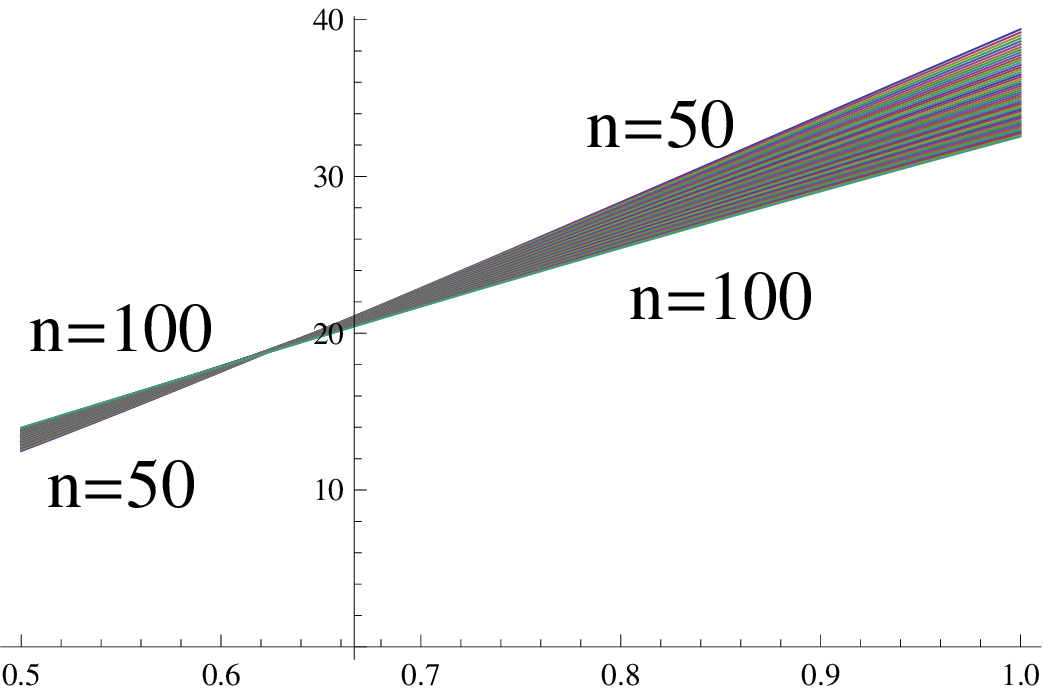}
		\hspace{0.2cm}
		\includegraphics[width=5cm,height=4cm]{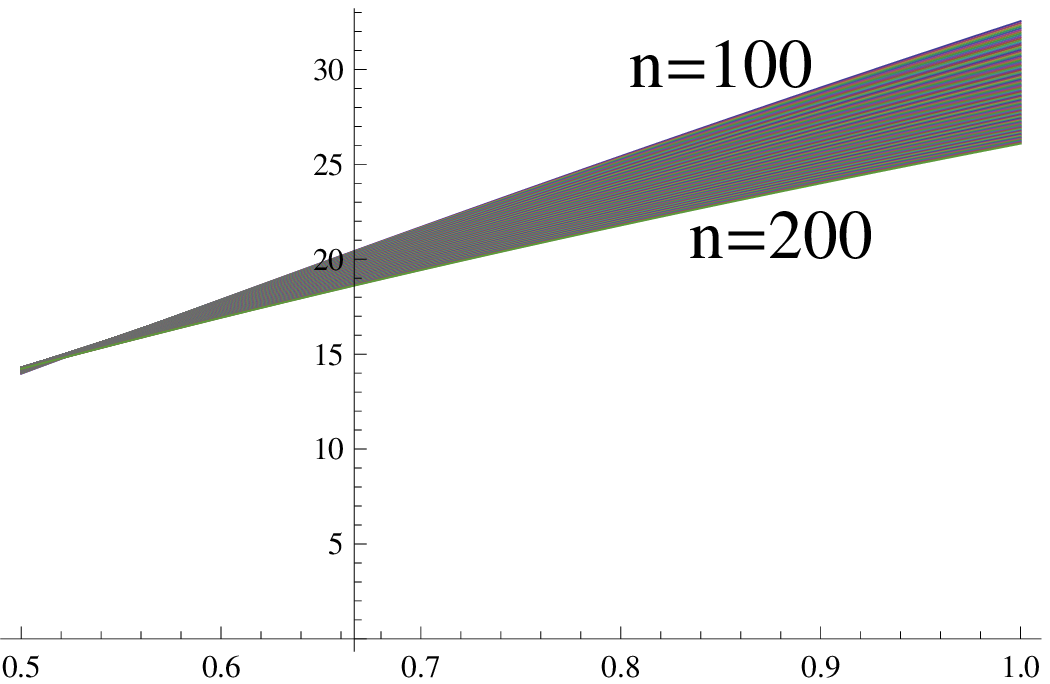}
	\caption{The curves $D_n(\tau)$ for $n\leq 50$ (left),
	$50\leq n\leq 100$ (center) and 
	$100\leq n\leq 200$ (right).}
\end{figure}

In our numerical experiments, we used the well known Shepp-Logan Phantom  with a slight modification
as shown on Figure 1: we have removed the thin layer around the head, which represents the skull, 
because it disappears too quickly by the smoothing procedure
and causes phases 1 and 3 to overlap, masking phase 2.
For more complicated functions $f$, such as most photographic pictures, phases $1$ and $3$ 
also tend to overlap for similar reasons. 
Indeed, these pictures often have details at the pixel scale, including electronic noise due to the captor. 
Since these details disappear early phase 3 begins immediately, therefore phase 2 cannot be observed.

\section{Extension to higher dimensions and higher order elements}

The results on approximation by anisotropic bidimensional piecewise linear
finite elements that we have exposed in \S 2 have been generalized in \cite{M}
to the case of elements of arbitrary order $m-1$
defined on partitions of $\Omega\subset \RR^d$ by simplices.
Here, the local error is defined as
$$
e_{m,T}(f):=\|f-I_{m,T}f\|_{L^P(T)},
$$
where $I_{m,T}$ denotes the local interpolation operator
on $\P_{m-1}$ for a $d$-dimensional simplex $T$. This operator
is defined by the condition
$$
I_{m,T}v(\gamma)=v(\gamma),
$$
for all points $\gamma\in T$ with barycentric coordinates in
the set $\{0, \frac 1 {m-1},\frac 2 {m-1},\cdots,1\}$. Then
one defines for any homogeneous polynomial $\bq\in \P_{m}$,
$$
K_{m,p}(\bq):=\inf_{|T|=1}e_{m,T}(\bq)_p.
$$
We refer to $K_{m,p}$ as the {\it shape function}. 
For piecewise linear elements in dimension two, i.e. $m=d=2$, 
we have observed that $K_p=K_{2,p}$ has the
special form given by \iref{Kpexp} which justifies the
introduction of the quantity $A_p(f)$. In a similar way, it can 
easily be proved that for piecewise linear elements in higher dimension,
i.e. $m=2$ and $d>2$, one has
$$
c_1 |{\rm det} (\bq)|^{1/d}\leq K_{2,p}(\bq) \leq  c_2 |{\rm det} (\bq)|^{1/d}.
$$
For piecewise quadratic elements in dimension two, i.e. $m=3$ and $d=2$, 
it was proved in \cite{M} that
$$
c_1 |{\rm disc} (\bq)|^{1/4}\leq K_{3,p}(\bq) \leq  c_2 |{\rm disc} (\bq)|^{1/4}.
$$
for any {\it homogeneous} polynomial $\bq\in \PP_{3}$, where
$$
{\rm disc}(a x^3+ b x^2 y+ c x y^2+ d y^3):= b^2 c^2 - 4 a c^3 - 4 b^3 d + 18 a b c d - 27 a^2 d^2.
$$
For other values of $m$ and $d$, equivalent expressions
of $K_{m,p}(\bq)$ in terms of polynomials in the coefficients of $\bq$
are available but of less simple form.

Defining the finite element interpolation error by an optimally adapted partition
$$
\sigma_N(f)_p:=\inf_{\#(\cT)\leq N}\|f-I_{m,\cT}f\|_{L^p},
$$
where $I_{m,\cT}$ is the global interpolation operator for the 
simplicial partition $\cT$, the following generalization of 
\iref{limsupest1} is proved in \cite{M}:
\be
\limsup_{N\to +\infty} N^{\frac m d}\sigma_N(f)_p\leq C_d\left\|K_{m,p}\(\frac {d^mf}{m!}\)\right\|_{L^\tau(\Omega)},
\;\; \frac 1 \tau:=\frac 1 p+\frac m d.
\label{limsupest2}
\ee
The constant $C_d$ is equal to $1$ when $d=2$ but larger
than $1$ when $d>2$ due to the impossibility of exactly tiling the
space with locally optimized simplices.
If $f$ is a $C^m$ function of $d$ variables, it is therefore natural to 
consider the quantity 
\be
A_{m,p}(f):=\|K_{m,p}(d^m f)\|_{L^\tau(\Omega)},
\;\; \frac 1 \tau:=\frac 1 p+\frac m d,
\ee
as a possible way to measuring anisotropic smoothness. For $d=2$ and piecewise
linear elements, we have seen in \S 2 that $A_{2,p}(f)$ is equivalent to the quantity $A_p(f)$.

Similarly to $A_p$ we are interested in the possible extension of $A_{m,p}$
to cartoon functions. We first introduce a generalisation 
of the notion of cartoon functions to higher piecewise smoothness $m$ and dimension $d$.
\begin{definition} 
\label{defcartoongen}
Let $m\geq 2$ and $d\geq 2$ be two integers. Let $\Omega\subset \R^d$ be an open set.
We say that a function $f$ defined
on $\Omega$ is a $C^m$ cartoon function if it is almost everywhere of the form 
$$
f=\sum_{1\leq i\leq k} f_i \Chi_{\Omega_i},
$$
where the $\Omega_i$ are disjoint open sets with piecewise $C^2$ boundary, no cusps (i.e. satisfying an interior and exterior cone condition), and such that $\overline \Omega = \cup_{i=1}^k \overline \Omega_i$.
Additionally, for each $1\leq i\leq k$, the function $f_i$ is assumed to be $C^m$ on $\overline \Omega_i$. 
\end{definition}

Let us consider a fixed cartoon function $f$ on a polyhedral domain $\Omega\subset\R^d$ (i.e. $\Omega$ is
such that $\overline \Omega$ is a closed polyhedron), and a decomposition 
$(\Omega_i)_{1\leq i\leq k}$ of $\Omega$ as in definition \ref{defcartoongen}. As before we define 
$ 
\Gamma := \bigcup_{1\leq i\leq k} \partial \Omega_i,
$
the union of the boundaries of the $\Omega_i$. Our assumptions on the sets $(\Omega_i)_{1\leq i\leq k}$ imply that $\Gamma$ 
is the union of a finite number of open hypersurfaces $(\Gamma_j)_{1\leq j\leq l}$, and of a set $\cP$ of dimension $d-2$.

As in \S 3, we now consider a sequence $f_N$ of piecewise linear approximations of $f$ on 
simplicial partitions $\cT_N$ of cardinality $N$. We distinguish two types of elements of $\cT_N$. A simplex $T\in \cT_N$ is called ``regular'' if $T\cap \Gamma=\emptyset$, and we denote the set of these simplices by $\cT_N^r$. Other simplices are called ``edgy'' and their set is denoted by $\cT_N^e$.
We can again split $\Omega$ according to
$$
\Omega:=(\cup_{T\in \cT_N^r}T) \cup (\cup_{T\in \cT_N^e}T)=\Omega_N^r \cup \Omega_N^e.
$$
Heuristically, if the partitions $\cT_N$ are built with approximation error minimisation in mind, the number of elements should be balanced between $\cT_N^r$ and $\cT_N^e$. The partition $\cT_N^r$ tends to cover most of the surface of $\Omega$, with simplices of diameter $\leq C  N^{-\frac 1 d}$, and $L^\infty$ approximation error $|f-f_N|\leq C N^{-\frac m d}$ (since we use $\cP_{m-1}$ elements).
On the other hand, since $f$ has discontinuities along $\Gamma$, the $L^\infty$ approximation error on $\cT_N^e$ does not tend to zero, and $\cT_N^e$ should thus be chosen so as to produce a thin layer around $\Gamma$. 
Let $h$ be the typical diameter of an element of $\cT_N^e$. Since the $\Gamma_j$ has
bounded curvature, this layer can be made of width $\cO(h^2)$
and therefore the layer around $\Gamma$ has volume bounded by
$h^2 \cH_{d-1}(\Gamma)$  up to a fixed multiplicative constant, 
where $\cH_{d-1}(\Gamma)$ is the $d-1$ dimensional
Hausdorff measure of $\Gamma$. On the other hand the minimal number of 
such elements of diameters $h$ needed to cover $\Gamma$ is 
bounded by $h^{1-d}\cH_{d-1}(\Gamma)$ up to a fixed multiplicative constant.
Eventually, we find that the layer around $\Gamma$ has volume bounded by $C N^{-\frac 2 {d-1}}$. 

Hence we have the following heuristic error estimate, for a well designed anisotropic partition:
\begin{eqnarray*}
\|f-f_N\|_{L^p(\Omega)} &\leq& \|f-f_N\|_{L^p(\Omega_N^r)} +\|f-f_N\|_{L^p(\Omega_N^e)}\\
&\leq& \|f-f_N\|_{L^\infty(\Omega_N^r)} |\Omega_N^r|^{\frac 1 p}+  \|f-f_N\|_{L^\infty(\Omega_N^e)} |\Omega_N^e|^{\frac 1 p}\\
& \leq & C (N^{-\frac m d}  + N^{\frac {-2} {p (d-1)}})
\end{eqnarray*}
This leads us to define a critical exponent
$$
p_c = p_c(m,d):=\frac{2d}{m (d-1)}.
$$
If one measures the error in $L^p$ norm with $p>p_c(m,d)$, then the contribution of the edge neighbourhood
$\Omega_N^e$ dominates, while if $p<p_c(m,d)$ it is negligible compared to the contribution
of the smooth region $\Omega_N^r$. 
For the critical exponent $p=p_c(m,d)$ the two terms have the same order, which makes the situation more interesting. Note in particular that $p_c(2,2) =2$, which is consistent with our previous analysis.

For $p\leq p_c(m,d)$, we obtain the approximation rate $N^{-m/d}$ which suggests that approximation
results such as \iref{limsupest2} should also apply to cartoon functions and that the 
quantity $A_{m,p}(f)$ should be finite for such functions. We again need to
use a regularization approach, for the same reasons as in \S 3.
For a given dimension $d$, we consider a radial nonnegative function $\vp$ 
of unit integral and supported in the unit ball of $\R^d$, and define for $\delta >0$
\be
\vp_\delta(z) := \frac 1 { \delta^d} \vp\left(\frac z \delta\right) \text{ and } f_\delta = f * \vp_\delta.
\ee

In order to define the quantities of involved in our conjecture, we need to introduce the second fundamental form of a hypersurface.
At any point $x\in \Gamma\sm \cP$ be denote by $\n(x)$ the unit normal to $\Gamma$. Note that since 
$\Gamma$ is piecewise $C^2$, the map $x\mapsto \n(x)$ is 
$C^1$ on $\Gamma\sm \cP$.
Furthermore, we define $T_x\Gamma := \n(x)^\perp$, the tangent space to $\Gamma$ at $x$.
In a neighbourhood of $x\in \Gamma\sm \cP$, the hypersurface $\Gamma$ admits a parametrization of the form
$$
u\in T_x\Gamma \mapsto x + u+ \lambda(u) \n(x) \in \Gamma \sm \cP,
$$
where $\lambda$ is a scalar valued $C^2$ function.
By definition, the second fundamental form of $\Gamma$ at the point $x$ is the quadratic form $\II_x$ associated
to $d^2\lambda(0)$ which is defined on $T_x \Gamma\times T_x \Gamma$.
Alternatively, for all $u,v\in T_x\Gamma$ we have $\II_x(u,v) := -\<\partial_u \n, v\>$.
The Gauss curvature $\kappa(x)$ is the determinant of $\II_x$, in any orthonormal basis of $T_x \Gamma$,
$$
\kappa(z):= \det \II_z.
$$
For example, in two space dimensions the tangent space $T_x\Gamma$ is one dimensional, and we simply have 
$\II_x(u,v) = \kappa(x) \<u,v\>$.
We also denote by $\sigma(z)\in \{0, \cdots, d-1\}$ the signature of the quadratic form $\II_z$, which is defined as the number of its positive eigenvalues. 

With $\tau$ such that $\frac 1 \tau :=\frac m d+\frac 1 p$, we 
define
$$
S_p(f) :=  \|K(d^m f )\|_{L^\tau(\Omega\sm \Gamma)}=A_p(f_{|\Omega\sm \Gamma}).
$$
We conjecture the following generalization to Theorem \ref{threg}.

\begin{conjecture}
There exists $d$ positive constants $C(k)$, $k\in \{0,\cdots, d-1\}$, that depend on  $\vp,p,m,d$, such that,
with 
$$
E_p(f):= \| C(\sigma) |\kappa|^{\frac m {2d}} [f] \|_{L^\tau(\Gamma)}=
\(\int_\Gamma \left |C(\sigma(z)) |\kappa(z)|^{\frac m {2d}} [f(z)] \right |^\tau dz\)^{\frac 1 \tau},
$$
we have
\label{thAmdp}
\begin{itemize}
\item If $p<p_c$ then 
$$
\lim_{\delta\to 0} A_{m,p}(f_\delta) = S_p(f).
$$ 
\item If $p=p_c$ then 
$$
\lim_{\delta \to 0} \left(A_{m,p}(f_\delta)\right)  = \left(S_p(f)^\tau +  E_p(f)^\tau\right)^{1/\tau}.
$$ 
\item If $p>p_c$ then
$$
\lim_{\delta \to 0}  \delta^{\frac 1 {p_c} - \frac 1 p} A_{m,p}(f_\delta) =  E_p(f).
$$ 
\end{itemize}
\end{conjecture}

In the remainder of this section, we give some 
arguments that justify this conjecture. Given a cartoon function $f$, 
we define the sets $\Omega_\delta$, $\Gamma_\delta^0$, 
$\Gamma_\delta$ and $\cP_\delta$ similarly to \S 3.
We need to perform an asymptotic analysis of the integral
\be
\label{tointegrate}
\int_\Omega K(d^m f_\delta)^\tau = \int_{\Omega_\delta} K(d^m f_\delta)^\tau+\int_{\cP_\delta} K(d^m f_\delta)^\tau
+\int_{\Gamma_\delta} K(d^m f_\delta)^\tau.
\ee
As in the proof of Theorem \ref{threg} the contribution of $\cP_\delta$ can be proved to be negligible 
compared with those of $\Omega_\delta$ and $\Gamma_\delta$ as $\delta \to 0$. The contribution of
$\Omega_\delta$ satisfies
$$
\lim_{\delta \to 0} \int_{\Omega_\delta} K(d^m f_\delta)^\tau =\int_{\Omega\sm \Gamma} K(d^m f)^\tau.
$$ 
The main difficulty lies again in the contribution of $\Gamma_\delta$. Let us define $\tau_c$ by 
$$
\frac 1 {\tau_c} := \frac m d+ \frac 1 {p_c}.
$$
The contribution of $\Gamma_\delta$ can be 
computed if one can establish an estimate generalizing \iref{estimK2d} according to 
\be
\label{estimK}
\left|  \delta^{\frac 1 {\tau_c}} K_{m,p}(d^m f_\delta(z)) - |[f](x)| |\kappa(x)|^{\frac m {2 d}} \Phi_{m,d,\sigma(x)}(u) \right|\leq \omega(\delta)
\ee
where $\omega(\delta)\to 0$ as $\delta\to 0$, $x\in \Gamma_\delta^0$, $u\in [-1,1]$, $z = x+\delta u \n(x)$,
and where the function $\Phi_{m,d,k} : [-1,1]\to \R$ only depends on $m,d,k$ and $\vp$. If
\iref{estimK} holds, we then easily derive that
$$
\lim_{\delta \to 0} \delta^{\frac {\tau}{\tau_c}-1}\int_{\Gamma_\delta} K(d^m f_\delta)^\tau=\int_\Gamma C(\sigma) |\kappa|^{\frac {\tau m} {2d}} |[f]|^\tau,
$$
with $C(k):=\int_{-1}^1|\Phi_{m,d,k}(u)|^\tau du$, which leads to the proof of the conjecture.

We do not have a general proof of \iref{estimK} for any $m$, $p$ and $d$. 
In the following, we justify its validity in
two particular cases for which the explicit expression of $K_{m,p}$ is known to us:
piecewise quadratic in two space dimensions ($d=2$ and $m=3$)
and piecewise linear  in any dimension ($m=2$).

\paragraph{Piecewise quadratic elements in two dimensions.}
For all $\delta>0$, $x\in \Gamma_\delta^0$ and $u\in [-1,1]$, let $\pi_{x,\delta,u}\in \P_3$
 be the homogeneous cubic polynomial on $\R^2$ corresponding to $d^3 f_\delta (x+\delta u \n(x))$.
Let also $\pi_{x,u}\in \P_3$ be the homogeneous cubic polynomial on $\R^2$ defined by 
$$
\pi_{x,u}(\lambda \n(x)+\mu \t(x)) = -\lambda(\Phi''(u) \lambda^2- 3\Phi'(u) \kappa(x) \mu^2)
$$
for all $(\lambda,\mu)\in\R^2$, where $\Phi$ is defined
by \iref{defchi}. For all $x\in \Gamma$, we denote by $M_{x,\delta}$ the (symmetric) linear map defined by 
$$
M_{x,\delta}\n(x) = \delta \n(x) \text{ and } M_{x,\delta} \t(x)= \sqrt \delta \t(x).
$$
Then, using a reasoning similar to the one used in the appendix of this paper, it can be proved that
\be
\label{limdiscD}
\|\pi_{x,\delta,u} \circ M_{x,\delta} -[f](x)\pi_{x,u}\|\leq \omega(\delta).
\ee
where $\lim_{\delta\to 0} \omega(\delta) = 0$ and the function $\omega$ depends only on $f$.
Furthermore, it is proved in \cite{M} that 
$$
K_{3,p}(\bq) = C \sqrt[4]{|\disc \bq|}
$$
where the positive constant $C$ depends on $p$ and the sign of $\disc \bq$. Combining this expression with \iref{limdiscD} proves the Estimate \iref{estimK} and thus the conjecture in the case $m=3$ and $d=2$.

\paragraph{Piecewise linear elements in any dimension.} We use the second 
fundamental form of the discontinuity set $\Gamma$ in order to evaluate
$d^m f_\delta$ on $\Gamma_\delta$.
Characteristic functions are one of the simplest types of cartoon functions. In that case, it is possible
to establish a simple relation between the second fundamental form of
the edge set and the second derivatives of $f$ in a distributional sense:
if $\Omega\subset \R^d$ is a bounded domain with smooth boundary $\Gamma$ and inward normal $\n$,
we then have for all $C^2$ test function $\psi$ and $u,v\in \R^d$ 
\be
-\int_\Omega \partial^2_{u,v}\psi = \int_\Gamma \<u,\n\>\<v,\n\> (\partial_\n \psi - \tr(\II) \psi) +\II'(u,v)\psi,
\label{d2charR2eqd}
\ee
where $\II'_x(u,v)$ is the second fundamental form $\II_x$ applied to the orthogonal projection of $u$ and $v$ on $T_x\Gamma$.
The proof of this formula (that generalizes \iref{d2charR2eq} which is proved in the appendix and
corresponds to the two-dimensional case $d=2$) is given further below. For all $x\in \Gamma$, 
we denote by $M_{x,\delta}$ the (symmetric) linear map defined by 
$$
M_{x,\delta}\n(x) = \delta \n(x) \text{ and } M_{x,\delta} \t = \sqrt \delta \t  
$$
for all $\t\in T_x\Gamma$.
For all $\delta>0$, $x\in \Gamma_\delta^0$ and $u\in [-1,1]$, let $\pi_{x,\delta,u}\in \P_2$ be the homogeneous quadratic polynomial on $\R^d$ corresponding to $d^2 f_\delta (x+\delta u \n(x))$.
Let also $\pi_{x,u}\in \P_2$ be the homogeneous quadratic polynomial on $\R^d$ defined by 
$$
\pi_{x,u}(\lambda \n(x)+\t) = \Phi'(u) \lambda^2- \Phi(u) \II_x(\t,\t)
$$
for all $\t\in T_x \Gamma$, where $\Phi(x):=\int_{\R^{d-1}}\vp(x,y)dy$.
Then, using formula \iref{d2charR2eqd} and a reasoning analogous to the one presented in the appendix, 
it can be proved that 
\be
\label{limQuadD}
\|\pi_{x,\delta,u} \circ M_{x,\delta} -[f](x)\pi_{x,u}\|\leq \omega(\delta).
\ee
where $\lim_{\delta\to 0} \omega(\delta) = 0$ and the function $\omega$ depends only on $f$.
Furthermore, it is proved in \cite{M} that 
$$
K_{2,p}(\bq) = C \sqrt[d]{|\det \bq|}
$$
where the positive constant $C$ depends on $d,p$ and the signature of $\bq$. Combining this expression with \iref{limQuadD} proves the estimate \iref{estimK} and thus the conjecture in the case $m=2$ in any dimension $d>1$.
\nl
\nl
{\bf Proof of \iref{d2charR2eqd}:}
Let $\proj_\Gamma$ be the orthogonal projection onto $\Gamma$, and for all $x\in \Gamma$ let $\proj_x$ be the orthogonal projection onto $T_x\Gamma$. 
We consider a vector $u\in \R^d$ and we define $\move : \Gamma\to\Gamma$ by $\move (x) := \proj_\Gamma(x+u)$. If $\|u\| \|\II\|_{L^\infty(\Gamma)}<1$, then $\move$ is smooth and its differential $d_x \move : T_x\Gamma\to T_{x'}\Gamma$, where 
$$
x' = \move(x),
$$
is given by the following formula
$$
d_x \move = (\Id - \<u, \n(x')\>\II_{x'})^{-1} \proj_{x'}. 
$$
The determinant of $d_x U$ (more precisely the determinant of the matrix of $d_x U$ in direct orthogonal bases of $T_x\Gamma$ and $T_{x'} \Gamma$) is 
$$
\det (d_x \move) = \det(\Id - \<u, \n(x')\>\II_{x'})^{-1}\<\n(x),\n(x')\> = 1+\<u, \n(x')\> \tr(\II_{x'})+ \|u\|\omega_1(u,x).
$$
where 
$\omega_1(u,x)$ tends uniformly to $0$ as $u\to 0$.
Furthermore, it is easy to show that 
$$
|\psi(x+u) - \psi(x') -  \<u, \n(x')\> \partial_{\n(x')} \psi(x')|\leq C\|u\|^2,
$$
and $\|n(x') - n(x) - \II_{x'}(\proj_{x'}(u))\|\leq \|u\|\omega_2(u)$, where $C$ and $\omega_2$ are independent of $x\in \Gamma$ and $\omega_2(u) \to 0$ as $u \to 0$. Combining these results, we obtain
$$
\int_\Gamma \psi(x+u)\<\n(x),v\> dx = \int_\Gamma \psi(x+u)\<\n(x),v\> \det (d_{x'} \move)^{-1} dx'\\
$$
$$
=\int_\Gamma \<\n(x'), v\> \psi(x')dx'+ \int_\Gamma \<\n(x'), v\> \<\n(x'),u\>(\partial_{\n(x')} \psi- \tr(\II_x')\psi(x'))
 +\<v,\II_x'(\proj_{x'}(u))\>\psi(x') dx' + \|u\|\omega_3(u),
$$
where $\omega_3(u)\to 0$ as $u \to 0$.
We conclude the proof of \iref{d2charR2eqd} using the formula
$$
- \int_\Omega \partial^2_{u,v}\psi = \lim_{h\to 0} h^{-1}\int_\Gamma (\psi(x+hu)-\psi(x))\<\n(x),v\> dx.
$$ 
\sq

\begin{remark}
Similarly to the results presented in \S4, there is an affine invariance property associated to $\kappa$: if $T$ is an affine transformation of $\R^d$ with linear part $L$, 
and if $f = \ti f\circ T$, $\ti \Gamma = T(\Gamma)$ and $\ti \kappa$ is the Gauss curvature of $\ti \Gamma$, then
one has for any $s \geq 0$, 
$$
(\det L)^{\frac{d-1}{d+1}}\int_\Gamma |C(\sigma) [f]|^s |\kappa|^{\frac 1 {d+1}} = \int_{\ti\Gamma} |C(\ti \sigma)[\ti f]|^s |\ti \kappa|^{\frac 1 {d+1}}.
$$
It follows from this observation that when $p=p_c$, the contribution of the edges is affine invariant in the sense
that 
$$
E_{p_c}(\ti f) = (\det L)^{\frac {d-1}{d+1}} E_{p_c}(f).
$$
Since one also has $A_{m,p}(\ti f) = (\det L)^{\frac {d-1}{d+1}} A_{m,p} (f)$ this comforts the conjecture.
Let us mention that the quantity $|\kappa|^{\frac 1 {d+1}}$ has been used in \cite{Ol} 
in order to define surface smoothing operators that are invariant under affine change of coordinates.
\end{remark}

\section{Conclusion}

In this paper we have investigated the quantity $A_p(f)$ which
governs the rate of approximation by anisotropic $\P_1$ finite elements as a way to
describe anisotropic smoothness of functions. This quantity is not a semi-norm
due to the presence of the non-linear quantity $\det(d^2f)$ 
and cannot be defined in a straightforward manner for general distributions.
We nevertheless have shown that this quantity can be defined
for cartoon images with geometrically smooth edges when $p\leq 2$.
A theoretical issue remains to give a satisfactory meaning to 
the full class of functions for which this quantity is finite.

From a more applied perspective, it could be interesting to
investigate the role of $A_p(f)$ in problems where anisotropic
features naturally arise: 

\begin{enumerate}
\item
Approximation of PDE's: in the case of one dimensional hyperbolic conservation laws, it was proved
in \cite{DLu} that despite the appearance of discontinuities
the solution has high order smoothness in Besov spaces
that govern the rate of adaptive approximation by piecewise polynomials.
A natural question is to ask wether similar results hold
in higher dimension, which corresponds to understanding
if $A_p(f)$ remains bounded despite the appearance of shocks.
\item
Image processing: as illustrated in \S 5, the quantity $A_p(f)$ can easily be discretized
and defined for pixelized images. It is therefore tempting to use
$A_2(f)$ in a similar way
as the total variation in \iref{minBV}, by solving a problem of the form
\be
\min_{g\in BV} \{A_2(g)\; ; \; \|Tg-h\|_{L^2}\leq \e\},
\label{minA2}
\ee
with the objective of promoting images with piecewise smooth edges. The main
difficulty is that $A_2$ is not a convex functional. One way to solve this
difficulty could be to reformulate \iref{minA2} in a Bayesian framework as the search of a maximum
of an a-posteriori probability distribution (MAP) as an estimator of $f$. In this framework, we
may instead search for a minimal mean-square error estimator (MMSE),
and this search can be implemented by stochastic algorithms which
does not require the convexity of $A_2$, see \cite{LMo}.
\end{enumerate}

\section*{Appendix: proof of the estimates \iref{estimnn2d}-\iref{estimnt2d}-\iref{estimtt2d}.}

It is known since
the work of Whitney on extension theorems 
(see in particular \cite{Whitney}) 
that for any open set $U\subset \R^d$, and any $g\in C^2(\overline U)$ there exists $\ti g \in C^2(\R^d)$ such that $\ti g_{|U} = g$.
It follows that for each $1\leq i \leq k$, there exists $\ti f_i\in C^2(\R^2)$, compactly supported, and such that $\ti f_{i|\Omega_i} = f_i$. 

Let $\Gamma_j$ be one of the pieces of $\Gamma$, between the domains $\Omega_k$ and $\Omega_l$, and let  $s=\ti f_k$ and $t = \ti f_l- \ti f_k$. Although the domains $\Omega_k$ and $\Omega_l$ are only piecewise smooth, there exists an open set $\Omega'$ with $C^2$ boundary such that for $\delta_0>0$ small enough
$$
f = s \Chi_{\Omega'} + t \ \text{ on } \ \bigcup_{0<\delta\leq \delta_0} (\Gamma_{j,\delta}+B_\delta),
$$
where $B_\delta$ is the ball of radius $\delta$ centered at $0$.
Note that $\Gamma_j\subset \Gamma' := \partial \Omega'$ and that $s = [f]$ on $\Gamma_j$.
In the following, the variables $x,z$ are always subject to the restriction
\be
\label{restrictXZ}
x\in \Gamma_{j,\delta}^0 \text{ and } z = U_\delta(x,u) = x+\delta u \n(x) \text{ where } 0<\delta\leq \delta_0 \text{ and } |u|\leq 1,
\ee
note that $z\in \Gamma_{j,\delta}$ and $\|x-z\|\leq \delta$. We therefore have 
$$
f_\delta(z) = \int_{\Omega'} s(\ti x) \vp_\delta(z-\ti x) d\ti x + t_\delta(z),
$$
where $t_\delta:=t * \vp_\delta$. The second derivatives of $t_\delta$ are uniformly bounded, and are therefore negligible in regard of all three estimates \iref{estimnn2d}, \iref{estimnt2d} and \iref{estimtt2d}, indeed
$$
\|d^2 t_\delta\|_{L^\infty} = \|(d^2 t)* \vp_\delta\|_{L^\infty} \leq \|d^2 t\|_{L^\infty} \|\vp_\delta\|_{L^1} = \|d^2 t\|_{L^\infty} \|\vp\|_{L^1} <\infty.
$$
We now define the $2\times 2$ symmetric matrices
$$
I(z,x) := \int_{\Omega'} (s(\ti x)-s(x)) \ d^2 \vp_\delta(z-\ti x) \ d\ti x \text{ and } J(z) := \int_{\Omega'} d^2 \vp_\delta(z-\ti x) d\ti x 
$$
so that 
\be
d^2 f_\delta(z) =  d^2 t_\delta + I(z,x) + [f](x)  J(z).
\label{d2fdelta}
\ee
We already know that the contribution of $d^2t_\delta$ is negligible. We now prove that
the same holds for the contribution of $I(z,x)$. Since $\vp_\delta(z-\ti x)$ is non-zero
only if $\|\ti x-z\| \leq \delta$ and therefore $\|\ti x-x\| \leq 2\delta$, we
can bound the norm of the matrix $I(z,x)$ by
\be
\label{d2f1}
\|I(z,x)\| \leq 2\delta \|ds\|_{L^\infty} \|d^2 \vp_\delta\|_{L^1} \leq 2\delta \|ds\|_{L^\infty} \|d^2\vp\|_{L^1} \delta^{-2} = C\delta^{-1}.
\ee
This proves that the contribution of $I(z,x)$ is negligible for the two estimates
\iref{estimnn2d} and \iref{estimnt2d}. In order to prove that it is
also negligible in the estimate \iref{estimtt2d},
we need a finer analysis of $\t(x)^\trans I(z,x) \t(x)$.
For this purpose we fix a unit vector $u$ and the pair $(x,z)$. We introduce
$$
\Lambda(\ti x) := (s(\ti x)-s(x)) \ \partial_u \vp_\delta(z-\ti x) + \partial_u s(\ti x) \ \vp_\delta(z-\ti x),
$$
so that by Leibniz rule
$$
(s(\ti x)-s(x)) \ \partial^2_{u,u} \vp_\delta(z-\ti x)= \partial^2_{u,u} s(\ti x) \ \vp_\delta(z-\ti x) - \partial_u \Lambda(\ti x).
$$
Therefore
\begin{eqnarray*}
u^\trans I(z,x)u &=& \int_{\Omega'} \left( \partial^2_{u,u} s(\ti x) \ \vp_\delta(z-\ti x) - \partial_u \Lambda(\ti x)\right) d\ti x, \\
&=&  \int_{\Omega'} \partial^2_{u,u} s(\ti x) \ \vp_\delta(z-\ti x) d\ti x - \int_{\Gamma'} \Lambda(\ti x) \<\n(\ti x), u\> d \ti x.
\end{eqnarray*}
The first integral clearly satisfies 
$$
\left |\int_{\Omega'} \partial^2_{u,u} s(\ti x) \ \vp_\delta(z-\ti x) d\ti x \right |\leq \|d^2 s\|_{L^\infty} \|\vp_\delta\|_{L^1},
$$
and is therefore bounded independently of $\delta$. We estimate the second integral 
for the special case $u = \t(x)$, remarking that  $|\<\n(\ti x), \t(x)\>|\leq C_1 \delta$ on the domain of integration.
Therefore
$$
\left| \int_{\Gamma'} \Lambda(\ti x) \<\n(\ti x), \t(x)\> d \ti x\right| \leq C_1 \delta |\Gamma'\cap B(z,\delta)| \|\Lambda\|_{L^\infty},
$$
where, slightly abusing notations, we denote by $|\Gamma'\cap B(z,\delta)|$ the length ($1$-dimensional Hausdorff measure) of the curve $\Gamma'\cap B(z,\delta)$.
Clearly $\Lambda(\ti x) = 0$ if $\|z-\ti x\|\geq \delta$. If $\|z-\ti x\|\leq \delta$ we have 
\be
\label{upperPsi}
|\Lambda(\ti x)|\leq (\|x-z\|+\|z-\ti x\|) \|ds\|_{L^\infty} \|d\vp\|_{L^\infty} \delta^{-3} + \|ds\|_{L^\infty} \|\vp\|_{L^\infty} \delta^{-2} \leq C_0\delta^{-2}
\ee
Since in addition $|\Gamma'\cap B(z,\delta)|\leq C_2 \delta$, we finally find that
$$
\left| \int_{\Gamma'} \Lambda(\ti x) \<\n(\ti x), \t(x)\> d \ti x\right| \leq C_0C_1C_2.
$$
We have therefore proved that 
 $$
|\t(x)^\trans I(z,x) \t(x)|\leq C,
 $$
 where the constant $C$ is independent of $\delta$, which shows that 
 the contribution of $I(z,x)$ is negligible in \iref{estimtt2d}.

We now analyze the contribution the quantity $[f]J(z)$ in \iref{d2fdelta}.
For this purpose, we use an expression of the second derivative of 
the characteristic function $\Chi_{\Omega'}$ of a smooth set $\Omega'$
in the distribution sense. 
We assume without loss of generality that $\Gamma'$ is parametrized 
in the trigonometric sense, and therefore that $\n$ is the inward normal to $\Omega$.
For all test function $\psi$, we have
$$
-\int_{\Omega'} \partial^2_{u,v} \psi = \int_{\Gamma'} \partial_u \psi \<v,\n\> =  \int_{\Gamma'} (\partial_\n \psi \<u,\n\> + \partial_\t\psi \<u,\t\>) \<v,\n\>
$$
and, by integration by parts, 
$$
\int_{\Gamma'} \partial_\t\psi \<u,\t\>\<v,\n\> = - \int_{\Gamma'} \psi \, (\<u,\kappa\n\>\<v,\n\>- \<u,\t\>\<v,\kappa\t\>). 
$$
Therefore, we have
\be
\label{d2charR2eq}
-\int_{\Omega'} \partial^2_{u,v} \psi = \int_{\Gamma'} \<u,\n\>\<v,\n\> (\partial_\n \psi-\kappa \psi) +\kappa \<u,\t\>\<v,\t\> \psi.
\ee
Applying this formula to $\psi(\ti x):=\vp_{\delta}(z-\ti x)$ we obtain
\be
-u^\trans J(z) v=  \int_{\Gamma'} \<u,\n(\ti x)\>\<v,\n(\ti x)\> (\partial_\n \vp_\delta(z-\ti x) - \kappa(\ti x)\vp_\delta(z-\ti x))  +\kappa(\ti x) \<u,\t(\ti x)\>\<v,\t(\ti x)\> \vp_\delta(z-\ti x) d\ti x
\label{uJv}
\ee
Since $\Gamma_j$ is $C^2$, there exists a constant $C_0$ such that for all $x_1,x_2\in \Gamma_j$, we have 
$$
|\<\t(x_1),\n(x_2)\>|\leq C_0 \|x_1- x_2\|,
$$
and 
$$
|1-\<\n(x_1),\n( x_2)\>| = |1-\<\t(x_1),\t( x_2)\>|\leq C_0 \|x_1- x_2\|^2.
$$
We finally remark that $ |\Gamma'\cap B(z,\delta)|\leq C_1 \delta$, and that 
$\|\vp_\delta\|_{L^\infty}\leq \|\vp\|_{L^\infty}\delta^{-2}$ and  $\|\partial_\n \vp_\delta\|_{L^\infty} \leq \|d\vp\|_{L^\infty} \delta^{-3}$. 

Taking the vectors $\t(x)$ or $\n(x)$ as possible values of $u$ and $v$ in \iref{uJv} and using the
above remarks, we obtain the estimates
\begin{eqnarray}
\label{threeineq1}
\left|\n(x)^\trans J(z)\n(x) + \int_{\Gamma'} \partial_\n \vp_\delta(z-\ti x) d \ti x\right| &\leq& C \delta^{-1},\\
\label{threeineq2}
|\t(x)^\trans J(z)\n(x)| &\leq& C \delta^{-1},\\
\label{threeineq3}
\left|\t(x)^\trans J(z)\t(x) + \int_{\Gamma'} \kappa(\ti x) \vp_\delta(z-\ti x)d \ti x\right| &\leq& C,
\end{eqnarray}
where the constant $C$ depends only on $f$. In view of \iref{d2fdelta}
We can immediately derive estimate
\iref{estimnt2d} from \iref{threeineq2}.

In order to derive the estimate \iref{estimtt2d} from \iref{threeineq3}, we first
introduce the modulus $\omega$ of continuity of $\kappa$ on $\Gamma_j$, 
$$
\omega(\delta) := \sup_{x_1, x_2\in \Gamma_j\sep \|x_1- x_2\|\leq \delta} |\kappa(x_1)-\kappa( x_2)|.
$$
Therefore
$$
\left|\t(x)^\trans J(z)\t(x) + \int_{\Gamma'} \kappa(\ti x) \vp_\delta(z-\ti x)d \ti x\right| 
\leq \left|\t(x)^\trans J(z)\t(x) + \kappa(x)\int_{\Gamma'}  \vp_\delta(z-\ti x)d \ti x\right| 
+ C\omega(\delta)\delta^{-1}.
$$
We now claim that 
\be
\left| \int_\Gamma \vp_\delta( z-\ti x) d\ti x - \delta^{-1} \Phi(u) \right| \leq C,
\label{phiPhi}
\ee
holds with $C$ independent of $\delta$ which implies the validity of  \iref{estimtt2d}. 
In order to prove \iref{phiPhi}, we use a local parametrization of $\Gamma'$: let 
$\lambda:\R\to \R$ be such that for $h$ small enough we have, $x+h\t(x)+\lambda(h)\n(x)\in \Gamma$. 
Note that we have $|\lambda(h)|\leq C_0 h^2$ and $|\lambda'(h)|\leq C_0h$ for $h$ small enough. 
Then for $\delta$ small enough,
\begin{eqnarray*}
\left| \int_{\Gamma'} \vp_\delta( z-\ti x) d\ti x - \delta^{-1} \Phi(u) \right|  &\leq&  \left|\int_{\R} \vp_\delta(h\t(x)+ (\delta u-\lambda(h))\n(x))\sqrt{1+\lambda'^2} dh - \int_{T_x\Gamma}  \vp_\delta(h\t(x)+ \delta u\n(x)) dh \right|\\
&\leq & C\delta(\|\vp_\delta\|_{L^\infty} (\sqrt{1+(C_0\delta)^2}-1)+  \|d\vp_\delta\|_{L^\infty} C_0 \delta^2 ) \leq C
\end{eqnarray*}
Finally, we can derive the estimate \iref{estimnn2d} from \iref{threeineq1} using
the inequality
\be
\left| \int_\Gamma \partial_\n \vp_\delta(z - \ti x) d\ti x + \delta^{-2} \Phi'(u)\right|  \leq C \delta^{-1}
\label{phiPhin}
\ee
which proof is very similar to the one of \iref{phiPhi}, using that $\|\n(x_1)-\n(x_2)\|\leq C \|x_1-x_2\|$ for all $x_1, x_2\in \Gamma_j$.

\noindent
\nl
\nl
Albert Cohen
\nl
UPMC Univ Paris 06, UMR 7598, Laboratoire Jacques-Louis Lions, F-75005, Paris, France
\nl
CNRS, UMR 7598, Laboratoire Jacques-Louis Lions, F-75005, Paris, France
\nl
cohen@ann.jussieu.fr

\noindent
\nl
\nl
Jean-Marie Mirebeau
\nl
UPMC Univ Paris 06, UMR 7598, Laboratoire Jacques-Louis Lions, F-75005, Paris, France
\nl
CNRS, UMR 7598, Laboratoire Jacques-Louis Lions, F-75005, Paris, France
\nl
mirebeau@ann.jussieu.fr

\end{document}